\documentclass[12pt,a4paper]{article} 
\usepackage{amsmath,pstricks,latexsym,theorem,epsfig} 
\usepackage{amsfonts,amssymb,latexsym,epsfig,amsmath}
\usepackage{amsmath,latexsym,theorem} 
\newtheorem{Theorem}{Theorem}[section] 
\newtheorem{Definition}{Definition}[section] 
\newtheorem{Proposition}{Proposition}[section] 
\newtheorem{Lemma}{Lemma}[section]

\theorembodyfont{\rmfamily}

\newcommand{\rec}[1]{{(\ref{#1})}} 
\newtheorem{Rm}{Remark}

\newcommand{\be}{\begin{equation}}
\newcommand{\ee}{\end{equation}}
\newcommand{\R}{\mathbb{R}}
\newcommand{\N}{\mathbb{N}}

\newcommand{\Z}{\mathbb{Z}}

\def\ti{\tilde}
\def\lf{\left}
\def\rg{\right}

\def\e{\varepsilon}
\def\ep{\varepsilon}
\def\ds{\displaystyle}

\def\Om{\Omega}
\def\om{\omega}

\def \N{I\!\!N} 

\def \R{I\!\!R}

\def \Z{Z\!\!\!Z} 
 
\def\11{1\!\!1}

 \def\res{\mathop{\hbox{\vrule height 7pt width .5pt 
depth 0pt\vrule height .5pt width 6pt depth 0pt}}\nolimits}
\def\ds{\displaystyle}

\def\rec#1{{(\ref{#1})}}

\newcommand{\ba}{\begin{array}} 
\newcommand{\ea}{\end{array}}

\begin{document} 
    \date{}
 \title{\sc Sub-criticality of non-local Schr\"odinger systems with antisymmetric potentials and applications to half-harmonic maps}
\author{ Francesca Da Lio\thanks{Department of Mathematics, ETH Z\"urich, R\"amistrasse 101, 8092 Z\"urich, Switzerland.} \thanks{Dipartimento di Matematica Pura ed Applicata, Universit\`a degli Studi di Padova. Via Trieste 63, 35121,Padova, Italy.} \and Tristan Riviere$^*$   }
\maketitle 
 
 \begin{abstract}
 We consider  nonlocal linear  Schr\"odinger-type critical systems    of the type
 \begin{equation}\label{eqabstr}
 \Delta^{1/4} v=\Omega\, v~~~\mbox{in  $\R\,.$} \
 \end{equation}
where $\Omega$ is antisymmetric potential in $L^2(\R,so(m))$, $v$ is a ${\R}^m$ valued map and $\Omega\, v$ denotes the matrix multiplication.  We show that every solution $v\in L^2(\R,\R^m)$ of \rec{eqabstr} is in fact in $L^p_{loc}(\R,\R^m)$, for every $2\le p<+\infty$, in other words, we prove that the system (\ref{eqabstr}) which is a-priori only critical in $L^2$ happens to have a subcritical behavior for antisymmetric potentials.
 As an application we obtain the $C^{0,\alpha}_{loc}$ regularity of weak  $1/2$-harmonic maps into $C^2$ compact manifold without boundary.
 
  \end{abstract}

  \tableofcontents
  \section{Introduction}
  In this paper we consider maps $v=(v_1,\cdots,v_m)\in L^2({\R},{\R}^m)$ solving a system of the form
 \begin{equation}
 \label{zz1}
\forall i=1\cdots m\quad\quad\quad\Delta^{1/4} v^i =\sum_{j=1}^m\Omega_j^i\, v^j\,,
\end{equation}
where $\Omega=(\Omega_i^j)_{i,j=1\cdots m}\in L^2({\R},so(m))$ is an $L^2$ maps from ${\R}$ into the space $so(m)$ of $m\times m$ antisymmetric matrices.  
The operator $\Delta^{1/4} $ on ${\R}$ is defined by means of the the Fourier transform as follows  $$\widehat{\Delta^{1/4} u}=|\xi|^{1/2}\hat u\,,$$
 (given a function f, $\hat f$ or ${\cal{F}}$ denotes the Fourier transform of $f$).\par
We will also simply denote such system in the following way
\[
\Delta^{1/4}v=\Omega\, v\,.
\]
 We remark that the system \rec{eqsch} is a-priori critical for $v\in L^2(\R)$. Indeed under the assumptions that $v,\Omega\in L^2$ we obtain that $\Delta^{1/4}v\in L^1$ and using
classical theory on singular integrals we deduce that $ v\in L^{2,\infty}_{loc}$, the weak-$L^2$ space, which has the same homogeneity of $L^2$.  Thus we are more or less back to the initial assumption which is a property that characterizes critical equations.\par

\medskip

In such critical situation it is a-priori not clear whether solutions have some additional regularity or whether weakly converging sequences of solutions tends
to another solution (stability of the equation under weak convergence)...etc. 

\medskip

In \cite{Riv} and \cite{Riv2} the second author proved the sub-criticality of local {\it a-priori} critical Sch\"odinger systems of the form
\begin{equation}
\label{zz2}
\forall i=1\cdots m\quad\quad\quad-\Delta u^i=\sum_{j=1}^m\Omega_j^i\cdot\nabla u^j\,,
\end{equation}
where $u=(u^1,\cdots,u^{m})\in W^{1,2}(D^2,{\R}^m)$ and $\Omega\in L^2(D^2,{\R}^2\times so(m))$ or of the form
\begin{equation}
\label{zz3}
\forall i=1\cdots m\quad\quad\quad-\Delta v^i=\sum_{j=1}^m\Omega_j^i\, v^j\,,
\end{equation}
where $v\in L^{n/(n-2)}(B^n,{\R}^m)$ and $\Omega\in L^{n/2}(B^n,{\R}^m)$. In each of these two situations the antisymmetry of $\Omega$
was responsible for the regularity of the solutions or for the stability of the system under weak convergence.

\medskip

Our first main result in this paper is to establish the sub-criticality of non-local Schr\"odinger systems of the form (\ref{zz1}). Precisely we prove
the following theorem which extends to a non-local setting the phenomena observed in \cite{Riv} and \cite{Riv2} for the above local systems.

  \begin{Theorem}\label{regschr}
Let $\Omega\in L^2(\R,so(m))$ and $v\in L^2(\R)$ be a weak solution of
\begin{equation}\label{eqsch}
\Delta^{1/4} v=\Omega\, v\,.
\end{equation}
  Then $v\in L^p_{loc}(\R)$ for every $1\le p<+\infty$.
\end{Theorem}

 As in the previous works the main technique for proving the above result is to proceed to a {\it change of gauge} by rewriting  the system
 after having multiplied $v$ by a well chosen rotation valued map $P\in H^{1/2}({\R},SO(m))$\footnote{$SO(m) $ is the space of $m\times m$ matrices $R$  satisying $R^tR=RR^t=Id$ and $det(R)=+1$} which is ''integrating'' $\Omega$ in an optimal
 way.  In \cite{Riv} the choice of $P$ for systems of the form (\ref{zz2}) was given by the geometrically relevant {\it Coulomb  Gauge} satisfying
 \be
 \label{zz4}
 div\lf[P^{-1}\nabla P+P^{-1}\Omega P\rg]=0\,.
 \ee
Here, like in \cite{Riv2}, an appropriate choice of the gauge $P$ satisfies a maybe less geometrically relevant equation which seems however
to be better adapted to the system (\ref{zz1}) :
\be
\label{zz5}
Asymm\lf[P^{-1}\Delta^{1/4} P\rg]=2^{-1}\lf[P^{-1}\Delta^{1/4} P  -\Delta^{1/4} P^{-1} P\rg]=\Omega \,.
\ee
The local existence of such $P$ is given by the following theorem.
\begin{Theorem}\label{P}
There exists $\varepsilon>0$ and $C>0$ such that for  
every $\Omega\in L^2(\R;so(m))$ satisying  $\int_{\R}|\Omega|^2 dx\le\varepsilon$,  there exists $P\in{ \dot {H}}^{1/2}(\R, SO(m))$  such that
\begin{equation}\label{cond}\begin{array}{ll}
(i)~~& \ds P^{-1}\Delta^{1/4} P  -\Delta^{1/4} P^{-1} P=2\,\Omega\,;
\\[5mm]
(ii)~~&\ds\int_{\R}|\Delta^{1/4} P|^2 dx\le C\int_{\R}|\Omega|^2 dx\,.
\end{array}
\end{equation}
\hfill $\Box$
\end{Theorem}
The proof of this theorem is established following an approach introduced by K.Uhlenbeck in \cite{Uhl}
while constructing {\it Coulomb Gauges} for $L^2$ curvatures
in 4 dimension. The construction does not provide the continuity of the map which to $\Om\in L^2$ assigns $P\in \dot{H}^{1/2}$. This illustrates
the difficulty of the proof of Theorem~\ref{P} which is not a direct consequence of an application of the local inversion theorem but requires more elaborated arguments.

\medskip
  
Thus if the $L^2$ norm of $\Omega$ is small, Theorem \ref{P}  gives a $P$ for which $w:=Pv$ satisfies
\begin{eqnarray}\label{Pveq}
\Delta^{1/4} w&=&-\left[P\Omega P^{-1}-\Delta^{1/4}P\, P^{-1}\rg]\,w+N(P,v)\nonumber\\
&=& -symm(\Delta^{1/4} P\, P^{-1})\ w+N(P,v)\,.
\end{eqnarray}
where $N$ is the bilinear operator defined as follows.
For an arbitrary integer $n$, for every $Q\in \dot{H}^{1/2}(\R^n, {\cal{M}}_{\ell\times m}(\R^n))$ $\ell\ge 0$\footnote{ ${\cal{M}}_{\ell\times m}(\R)$ denotes, as usual, the space of $\ell\times m$ real matrices.} and $v\in L^2(\R^n,\R^m)$,   $N$ is  given by
\begin{equation}\label{opTintr}
N(Q,v):=\Delta^{1/4}(Q\, v)-Q\Delta^{1/4} v+\Delta^{1/4} Q\ v\,.
\end{equation}
 One of the key result used in \cite{DLR} establishes that, under the above assumptions on $Q\in H^{1/2}({\R}^n,M_m({\R}))$ and $v\in L^2({\R}^n,{\R}^m)$,  $N(Q,v)$ is more regular than each of it's three
 generating terms respectively $\Delta^{1/4}(Q\ v)$, $Q\Delta^{1/4} v$ and $\Delta^{1/4} Q\, v$ \footnote{The last one for example being only  a-priori in $L^1$.}.  We proved that $N(Q,v)$ is
 in fact in $H^{-1/2}({\R},{\R}^m)$. Such a result in \cite{DLR} was called a 3-commutator estimate (see Theorem \ref{commDLR}).
 \medskip
 
 In the paper \cite{DLR2}  we are improve the gain of regularity by compensation
 obtained in \cite{DLR}.
 In order to make it more precise we recall the definition of the Hardy space
${\mathcal H}^1({\R}^n)$ which is the space of $L^1$
functions $f$ on ${\R}^n$satisfying 
\[
\int_{{\R}^n}\sup_{t\in {\R}}|\phi_t\ast f|(x)\ dx<+\infty\quad ,
\]
where $\phi_t(x):=t^{-n}\ \phi(t^{-1}x)$ and where $\phi$ is some function in the Schwartz space ${\mathcal S}({\R}^n)$ satisfying $\int_{{\R}^n}\phi(x)\ dx=1$.
\footnote{For more properties on the Hardy space ${\mathcal H}^1$ we refer to \cite{Gra1} and \cite{Gra2}.} 
\begin{Lemma}
\label{comm1new}
There exists a constant $C>0$ such that, for any $Q\in \dot{H}^{1/2}({\R}^n,M_m({\R}))$ and $v\in L^2({\R}^n,{\R}^m)$,  $N(Q,v)=\Delta^{1/4}(Q\, v)-Q\Delta^{1/4} v+\Delta^{1/4} Q\ v$ is in ${\mathcal H}^1({\R}^n)$
 and the following estimate holds
\begin{equation}\label{commest2}
\|N(Q,v)\|_{{\cal{H}}^{1}}\le C\, \|Q\|_{\dot{H}^{1/2}}\ \|v\|_{L^2({\R})}\,. 
\end{equation}
\end{Lemma}
Thus in equation (\ref{Pveq}) the last term in the r.h.s happens to be slightly more regular. It remains to deal with the first term in this r.h.s. : $-symm(\Delta^{1/4} P P^{-1})\ w$.
A-priori $symm(\Delta^{1/4} P P^{-1})=2^{-1}[\Delta^{1/4} P\, P^{-1}+P\,\Delta^{1/4} P^{-1}]$ is only in $L^2$ but here again we are going to take advantage of  a gain of regularity
 due to a compensation. Though, individually each of the terms $\Delta^{1/4} P\, P^{-1}$ and it's transposed $P\,\Delta^{1/4} P^{-1}$ are only in $L^2$, the sum happens to belong to the ''slightly'' smaller space $L^{2,1}$ defined as follows:
$L^{2,1}(\R)$ is the Lorentz space of measurable functions satisfying
$$
\int_{\R_+}t^{-1/2}f^*(t)dt<+\infty \,,$$
where $f^*$ is the decreasing rearrangement of $|f|\,.$

\medskip

The fact that $symm(\Delta^{1/4} P\, P^{-1})$ belongs to $L^{2,1}({\R})$ comes from the combination of the following lemma according to which
 $\Delta^{1/4}(symm(\Delta^{1/4} P\, P^{-1}))\in {\mathcal H}^1({\R})$ and the sharp Sobolev embedding \footnote{The fact that  $v\in  {\cal{H}}^1$ implies $\Delta^{-1/4} v\in L^{2,1}$  is deduced by duality from the fact that $\Delta^{1/4}v\in L^{2,\infty}$ implies that $v\in BMO({\R})$ - This last embedding has been proved by Adams in \cite{AD}} which says that $f\in {\mathcal H}^1({\R})$
implies that $\Delta^{-1/4}f\in L^{2,1}$. Precisely we have
\begin{Lemma}
\label{lm-I.2}
Let $P\in H^{1/2}({\R},SO(m))$ then $\Delta^{1/4}[symm(\Delta^{1/4} P\ P^{-1})] $
is in the Hardy space  ${\cal{H}}^1({\R})$ and the following estimates hold
\[
\|\Delta^{1/4}[\Delta^{1/4} P\ P^{-1}+P\ \Delta^{1/4} P^{-1}]\|_{{\mathcal H}^1}\le C\|P\|_{H^{1/2}}^2
\]
where $C>0$ is a constant independent of $P$. This implies in particular that
\[
\|symm(\Delta^{1/4} P\ P^{-1})]\|_{L^{2,1}}\le C\|P\|_{H^{1/2}}^2\,.
\]

\end{Lemma}
The proof of this lemma is a consequence of the {\it 3-commutator estimates} in \cite{DLR} (see Theorem \ref{commDLR3bis} below).

\begin{Rm}
\label{rm-I.1}
The fact that, for rotation valued maps $P\in W^{2,n/2}({\R}^n,SO(m))$ ($n>2$),  $symm(\Delta P\ P^{-1})$ happens to be more regular than $Asymm(\Delta P\ P^{-1})$ was also one of the key point
in \cite{Riv2}.
\end{Rm}

\medskip

As we explain in Section 3 Theorem~\ref{regschr} is a consequence of this special choice of $P$ for which the new r.h.s. in the gauge transformed equation~\rec{Pveq} is slightly
more regular due to lemmas~\ref{comm1new} and lemma~\ref{lm-I.2}.
More precisely this {\em gain of regularity} in the right of equation \rec{Pveq} combined with suitable localization arguments  permit to obtain the following local Morrey type estimate for $Pv$ and thus for $v$, being $P$ bounded,
\begin{equation}
\label{lqcondintr}
\sup_{ {{\displaystyle{\mathop{\scriptstyle{x_0\in B(0,\rho)}}_{0<r<\rho/4}}}}} r^{-\beta}\ \int_{B(x_0,r)}|\Delta^{1/4} v| dx\le C
\end{equation}
for $\rho$ small enough and $0<\beta<1/2$ independent on $x_0$.  Theorem 5.1 in \cite{AD} yields that $v\in L^q_{loc}(\R)$ for some $q>2\,.$ \footnote{ In a  paper in preparation \cite{DLR2} we show that the solutions of \rec{eqsch} are actually in $L^{\infty}_{loc}(\R)\,.$}

\medskip

Our study of the linear systems has been originally motivated by the following non-linear problem.

In the joint paper \cite{DLR} we proved the $C^{0,\alpha}_{loc}$ regularity of weak $1/2$-harmonic maps into a sphere $S^{m-1}$. The second aim of the present paper 
  is to extend this result to weak $1/2$-harmonic maps with values in a $k$ dimensional sub-manifold ${\cal{N}}$, which is supposed at least $C^2$,
  compact and without boundary. We recall that $1/2$-harmonic maps are functions $u$ in the space
  $
 \dot{H}^{1/2}(\R,{\cal{N}})=\{u\in \dot{H}^{1/2}(\R,\R^m):~~u(x)\in {\cal{N}}, {\rm a.e},\}\,,
 $
which are critical points for perturbation of the type $\Pi^{N}_ {\cal{N}}(u+t\varphi)$, ($\varphi\in C^{\infty}$ and
$\Pi^{N}_ {\cal{N}}$ is the normal projection on $ {\cal{N}}$)
 of the functional
\begin{equation}\label{lagr}
{\cal{L}}(u)=\int_{\R}|\Delta^{1/4} u(x)|^2 dx\,,
\end{equation}
  (see Definition 1.1 in \cite{DLR})\,.
 The Euler Lagrange equation associated to this non linear problem can be written as follows :
 
 \be
 \label{zz6}
 \Delta^{1/2}u\wedge \nu(u)=0\quad\quad\mbox{ in }\cal{D}'({\R})\,,\
 \ee
where $\nu(z)$ is the Gauss Maps at $z\in\cal{N}$ taking values into the grassmannian $\ti{G}r_{m-k}({\R}^m)$ of oriented $m-k$ planes in ${\R}^m$ which
 is given by the oriented normal $m-k-$plane to $T_z\cal{N}\,.$ \footnote{Since we are assuming that $\cal{N}$ is $C^2$, $\nu$ is a $C^1$ map on $\cal{N}$
 and the paracomposition gives that $\nu(u)$ is in $\dot{H}^{1/2}({\R},\wedge^{m-k}{\R}^m)$ hence, since $ \Delta^{1/2}u$ is {\it a-priori} in $\dot{H}^{-1/2}$ the product
 $\Delta^{1/2}u\wedge \nu(u)$ makes sense in $\cal{D}'$ using the duality $\dot{H}^{1/2}-\dot{H}^{-1/2}$} 
 
 \medskip
 
The  Euler Lagrange equation in the form (\ref{zz6}) is hiding fundamental properties of this equation such as it's elliptic nature...etc and is difficult to use directly 
for solving problems related to regularity and compactness. One of the first task is then to rewrite it in a form that will make some of it's
analysis features more apparent. This is the purpose of the next proposition. Before to state it we need some additional notations
 
 \medskip
 
Denote by $P^T(z)$ and $P^N(z)$  the projections respectively to the tangent space $T_z\cal{N}$ and to the normal space $N_z\cal{N}$ to  ${\cal{N}}$ at $z\in \cal{N}$ . For $u\in \dot{H}^{1/2}({\R},\cal{N})$ we simply denote by $P^T$ and $P^N$ the compositions $P^T\circ u$ and $P^N\circ u$. In Section 
\ref{harmonic}  we establish that , under the assumption $\cal{N}$ to be $C^2$, $P^T\circ u$ as well as $P^N\circ u$ are matrix valued maps in $\dot{H}^{1/2}({\R},M_m({\R}))$. 

\medskip

A useful formulation of the 1/2-harmonic map equation is given by the following result  

 \begin{Proposition}\label{EulEq}
Let $u\in\dot{H}^{1/2}(\R,{\cal{N}} )$ be a weak $1/2$-harmonic map. Then the following equation holds
\be\label{Euler}
\Delta^{1/4}v= \ti{\Om}_1+\ti{\Om}_2 \,v+\Om\,v \,,
\ee
where $v\in L^2({\R},{\R}^{2m})$ is given by
\[
v:=\left(\begin{array}{l}
P^T\Delta^{1/4} u\\
{\cal{R}}P^N\Delta^{1/4} u\end{array}\right)
\]
and where ${\cal{R}}$ is the Fourier multiplier of symbol $\sigma(\xi)=i\frac{\xi}{|{\xi}|}$. 

\noindent $\Om\in L^2({\R},so(2m))$ is given by
$$
\Omega=2\left(\begin{array}{cc}
- \omega &  \omega_{{\cal{R}}}\\[5mm]
 \omega_{{\cal{R}}}  & - {\cal{R}} \omega_{{\cal{R}}} \end{array}\right)$$the maps $\omega$ and $\omega_{{\cal{R}}}$ are in $L^2({\R},so(m))$ and given respectively by
$$
\omega=\frac{\Delta^{1/4} P^T P^T- P^T\Delta^{1/4} P^T}{2}\,,
$$
and
$$
\omega_{{\cal{R}}}=\frac{({\cal{R}}\Delta^{1/4} P^T) P^T- P^T({\cal{R}}\Delta^{1/4} P^T)}{2}\,.
$$
Finally the maps $\ti{\Om}_1$ and $\ti{\Om}_2$ are respectively in $L^{2,1}({\R},M_{2m}({\R}))$ and in ${\cal{H}}^1({\R},{\R}^{2m})$. 
\end{Proposition}\par

The explicit formulations of $\ti\Omega_1$ and $\ti{\Om}_2$  are given in  Section \ref{harmonic}.  
The control on $\ti\Omega_1$ and $\ti{\Om}_2$ is a consequence of regularity by compensation results on some operators that we now introduce.
\par
For every  $Q,v\in L^2(\R^n)$ we define the operator $ F$   by
\begin{equation}\label{F}
F(Q,v):={\mathcal{R}}(Q){\mathcal{R}}(v)-Qv\,.
\end{equation}

 In \cite{CRW} it is shown that  $F(Q,v)\in 
H^{-1/2}(R)$ and
\begin{equation}\label{estFinfty}
\|F(Q,v)\|_{ H^{-1/2}(R)}\le C \|Q\|_{L^2(\R)}\|v\|_{L^{2}(\R)}\,.
\end{equation}
 By  a suitable estimate on the dual operator of $F$  (Lemma \ref{Romega})   we deduce the  following sharper  estimate 
\begin{equation}\label{estFinftybis}
\|F(Q,v)\|_{ H^{-1/2}(R)}\le C \|Q\|_{L^2(\R)}\|v\|_{L^{2,\infty}(\R)}\,.
\end{equation}

 Next we recall some {\em commutator estimates}  we obtained in \cite{DLR}.

 %%%%%%%%%%%%%%%%%%%%%
  \begin{Theorem}\label{commDLR} Let $n\in {\N}^\ast$ and
 let $u\in BMO(\R^n)$, $Q\in \dot{H}^{1/2}({\R}^n,{\cal{M}}_{\ell\times m}(\R))$ . Denote
 \[
 T(Q,u):=\Delta^{1/4}(Q\Delta^{1/4}u)-Q\Delta^{1/2} u+\Delta^{1/4} u\Delta^{1/4} Q\quad ,
 \]
 then $T(Q,u)\in {H}^{-1/2}(\R^n)$ and there exists $C>0$, depending only on $n,$ such that
\begin{equation}
\label{zz7bis}
\|T(Q,u)\|_{{{H}}^{-1/2}(\R^n)}\le C\ \|Q\|_{\dot{H}^{1/2}(\R^n)}\ \|u\|_{BMO(\R^n)}\,.
\end{equation}
\hfill$\Box$
\end{Theorem}
  \begin{Theorem}\label{commDLR2} Let $n\in{\N}^\ast$ and 
 let $u\in BMO(\R^n)$, $Q\in \dot{H}^{1/2}({\R}^n,{\cal{M}}_{\ell\times m}(\R))$ . Denote
 \[
 S(Q,u):=\Delta^{1/4}[Q\Delta^{1/4} u]-{\cal{R}}  (Q\nabla u)+{\cal{R}}(\Delta ^{1/4} Q{\cal{R}}\Delta ^{1/4} u)
 \]
 where ${\cal{R}}$ is  the Fourier multiplier of symbol $m(\xi)=i\frac{\xi}{|{\xi}|}\,$. Then $S(Q,u)\in{H}^{-1/2}(\R^n)$ and there exists $C$ depending only on $n$ such that
\begin{equation}
\label{zz8bis}
\|S(Q,u)\|_{{{H}}^{-1/2}(\R^n)}\le C\ \|Q\|_{\dot{H}^{1/2}(\R^n)}\|u\|_{BMO(\R^n)}\,.
\end{equation}
 \hfill $\Box$
\end{Theorem}
%%%%%%%%%%%%%%%%
As it is observed in \cite{DLR} Theorems \ref{commDLR} and \ref{commDLR2} are consequences respectively  of
   the following results which are their ``dual versions"\,.
\begin{Theorem}\label{commDLR3bis}
Let $u,Q\in \dot H^{1/2}(\R^n)$, denote
$$
R(Q,u)=\Delta^{1/4}(Q\Delta^{1/4} u)-\Delta^{1/2}(Qu)+\Delta^{1/4}((\Delta^{1/4} Q) u)\,.$$
then $R(Q,u)\in {\cal{H}}^{1}(\R^n)$ and
\begin{equation}\label{zz9bis}
\|R(Q,u)\|_{{\cal{H}}^{1}(\R^n)}\le C\|Q\|_{\dot{H}^{1/2}(\R^n)}\|u\|_{\dot{H}^{1/2}(\R^n)}\,.\end{equation}
\end{Theorem}
\begin{Theorem}\label{commDLR4bis}
Let $u,Q\in \dot H^{1/2}(\R^n)$, denote $$
\tilde S(Q,u)=\Delta^{1/4}(Q\Delta^{1/4} u)-\nabla(Q {\cal{R}}u)+{\cal{R}}\Delta^{1/4}(\Delta^{1/4} Q {\cal{R}}u)\,.$$
 Then $\tilde S(Q,u)\in {\cal{H}}^1(\R^n)$ and 
\begin{equation}\label{zz10bis}
\|\tilde S(Q,u)\|_{{\cal{H}}^{1}(\R^n)}\le C\|Q\|_{\dot{H}^{1/2}(\R^n)}\|u\|_{\dot{H}^{1/2}(\R^n)}\,.\end{equation}
\hfill $\Box$
\end{Theorem}
  Since the operators $R$ and $\tilde S$ are the duals respectively  of $T$ and $S$, by combining
  Theorems \ref{commDLR} and \ref{commDLR3bis} and Theorems \ref{commDLR2} and 
  \ref{commDLR4bis}
  one gets the followings sharper estimates for $T$ and $S$:
  \begin{equation}
\label{zz7tris}
\|T(Q,u)\|_{{{H}}^{-1/2}(\R^n)}\le C\ \|Q\|_{\dot{H}^{1/2}(\R^n)}\|\Delta^{1/4}u\|_{L^{2,\infty}(\R^n)}\,;
\end{equation}
\begin{equation}
\label{zz8tris}
\|S(Q,u)\|_{{{H}}^{-1/2}(\R^n)}\le C\ \|Q\|_{\dot{H}^{1/2}(\R^n)} \|\Delta^{1/4}u\|_{L^{2,\infty}(\R^n)}\,.
\end{equation}

An adaptation of theorem~\ref{regschr} to the Euler Lagrange equation of the 1/2-Energy written in the form (\ref{Euler}) leads to the following theorem
which is the second main result of the present paper.  
 
 \begin{Theorem}\label{reghm}
Let  ${\cal{N}}$ be a closed $C^2$ submanifold of ${\R}^m$. Let  $u\in \dot{H}^{1/2}(\R,{\cal{N}})$ be a weak $1/2-$harmonic map into ${\cal{N}}$, then  $u\in C^{0,\alpha}_{loc}(\R,{\cal{N}})\,.$ \hfill$\Box$
\end{Theorem}
Finally a classical elliptic type bootstrap argument leads to the following result  
(see \cite{DLR2} for the details of this argument).
\begin{Theorem}
\label{th-I.3}
Let ${\cal{N}}$ be a smooth closed submanifold of ${\R}^m$. Let $u$ be a weak $1/2$-harmonic map in $\dot{H}^{1/2}({\R},{\cal{N}}))$, then $u$ is $C^\infty$ .
\hfill $\Box$
\end{Theorem}

The regularity of critical points of non-local functionals  has been recently investigated by Moser \cite{Moser}. In this work critical points to the functional
that assigns to any $u\in \dot{H}^{1/2}({\R},{\cal{N}})$ the minimal Dirichlet energy among all possible extensions \underbar{ in $\cal{N}$} are considered,
while in the present paper the classical $\dot{H}^{1/2}$ Lagrangian corresponds to the minimal Dirichlet energy among all possible extensions \underbar{ in ${\R}^m$}.
Hence the approach in \cite{Moser} consists in working with an \underbar{intrinsic} version of $H^{1/2}-$energy while we are considering here an \underbar{extrinsic} one.
The drawback of considering the \underbar{intrinsic} energy is that the Euler Lagrange equation is almost impossible to write explicitly and is then \underbar{implicit} while in the present case it has the
\underbar{explicit} form (\ref{zz6}). However the intrinsic version of the $1/2-$harmonic map is more closely related to the existing regularity theory of Dirichlet Energy minimizing maps into $\cal{N}$.

\medskip

The paper is organized as follows.

\begin{itemize}

\item[-] In Section \ref{regularitySchr} we prove Theorem \ref{regschr} .\par
\item[-] In Section \ref{constrP} we prove Theorem  \ref{P}\,.\par
\item[-] In Section \ref{harmonic} we  derive the Euler-Lagrange equation (\ref{Euler}) associated to the Lagrangian 
\rec{lagr} and we prove Theorem \ref{reghm}\,.

\item[-] In Appendix \ref{LocEnergy} we prove some $L-$energy decrease control  for   solutions to 
  linear non-local Schr\"ondiger type systems .    
 \item[-] In Appendix \ref{commutators} we provide commutator estimates that are crucial for the construction of the gauge $P$.  
\end{itemize}

\section{Preliminaries: function spaces and the fractional Laplacian}\label{defnot}

In this Section we introduce some  notations and definitions we are going to use in the sequel.
\par

For $n\ge 1$,  we denote respectively by ${\cal{S}}(\R^n)$ and  ${\cal{S}}^{\prime}(\R^n)$ the spaces of Schwartz functions and tempered distributions. Moreover given a function $v$ we will denote either by
  $\hat v$ or by  ${\mathcal{F}}[v]$ the Fourier Transform of  $v$ :
  $$\hat v(\xi)={\mathcal{F}}[v](\xi)=\int_{\R^n}v(x)e^{-i \langle \xi, x\rangle }\,dx\,.$$
  Throughout the paper we use the convention that $x,y$ denote variables in the space and
  $\xi$ the variable in the phase\,.
  
  We recall the definition of fractional Sobolev space (see for instance \cite{T3}).\par
  \begin{Definition}\label{fracsob} 
  For a real $s\ge 0$, 
  $$H^{s}(\R^n)=\lf\{v\in L^2(\R^n):~~|\xi|^s{\cal{F}}[v]\in L^2(\R^n)\rg\}\quad$$
  For a real $s<0$,
   $$H^{s}(\R^n)=\lf\{v\in {\cal{S}}^{\prime} (\R^n):~~(1+|\xi|^2)^s{\cal{F}}[v]\in L^2(\R^n)\rg\}\,.$$
   \hfill $\Box$
   \end{Definition}
   It is known that $H^{-s}(\R^n)$ is the dual of $H^{s}(\R^n)\,.$\par
   For $0<s<1$ another classical characterization of $H^{s}(\R^n)$which does not make use the Fourier transform is the following, (see for instance \cite{T3}).
   \begin{Lemma}
   For $0<s<1$, $u\in H^{s}(\R^n)$ is equivalent to $u\in L^2(\R^n)$ and
   $$
 \left(\int_{\R^n}\int_{\R^n} \left(\frac{(u(x)-u(y))^2}{|x-y|^{n+2s}}\right)dx dy\right)^{1/2}<+\infty\,.$$
 \hfill $\Box$\end{Lemma}
  For $s>0$ we set
  $$
  \|u\|_{H^s(\R^n)}=\|u\|_{L^2(\R^n)}+\\|\xi|^s{\cal{F}}[v]\|_{L^2(\R^n)}\,,
  $$
  and
  $$
  \|u\|_{\dot{H}^s(\R^n)}=\\|\xi|^s{\cal{F}}[v]\|_{L^2(\R^n)}\,.
  $$
  For an open set $\Omega\subset\R^n$,   $H^s(\Omega)$ is the space of the restrictions of functions from $H^{s}(\R^n)$ and
 $$
 \|u\|_{\dot H^{s}(\Omega)}=\inf\lf\{\|U\|_{\dot H^{s}(\R^n)},~~U=u ~\mbox{on}~\Omega\rg\}
 $$
 In the case  $0<s<1$ then 
  $u\in H^{s}(\Omega)$ if and only if 
 $u\in L^2(\Omega)$ and  
 $$
 \left(\int_\Omega\int_\Omega \left(\frac{(u(x)-u(y))^2}{|x-y|^{n+2s}}\right)dx dy\right)^{1/2}<+\infty\,.
 $$
 Moreover
 $$
  \|u\|_{\dot H^{s}(\Omega)}\simeq \left(\int_\Omega\int_\Omega \left(\frac{(u(x)-u(y))^2}{|x-y|^{n+2s}}\right)dx dy\right)^{1/2}<+\infty\,,
  $$
  see for instance \cite{T3}\,.
  
  Finally for a submanifold  ${\cal{N}}$ of $\R^m$  we can define   
 $$H^{s}({\R}^n,{\cal{N}})=\{u\in H^{s}({\R}^n,\R^m):~~u(x)\in {\cal{N}}, {\rm a.e.}\}\,.
 $$
 Given $q>2$ we also set
 $$
 W^{s,q}(\R^n):=\{v\in L^q({\R}^n):~~|\xi|^s{\cal{F}}[v]\in L^q(\R^n)\}\,.
 $$
 
 \medskip

 We shall make use of the Littlewood-Paley  dyadic  decomposition of unity that we recall here. Such a decomposition can be obtained as follows \,.
 Let $\phi(\xi)$ be a radial Schwartz function supported in $\{\xi\in{\R}^n:~|\xi|\le 2\}$, which is
 equal to $1$ in $\{\xi\in{\R}^n: ~|\xi|\le 1\}$\,.
 Let $\psi(\xi)$ be the function given by
 $$
 \psi(\xi):=\phi(\xi)-\phi(2\xi)\,.
 $$
 $\psi$ is then a ''bump function'' supported in the annulus $\{\xi\in{\R}^n:~1/2\le |\xi|\le 2\}\,.$
 
 \medskip
 
 Let $\psi_0=\phi$, $\psi_j(\xi)=\psi(2^{-j}\xi)$ for $j\ne 0 \,.$ The functions $\psi_j$, for $j\in\Z$, are supported in  $\{\xi\in{\R}^n:~2^{j-1}\le |\xi|\le 2^{j+1}\}\,$
and they realize  a dyadic decomposition of  the unity :
  $$
  \sum_{j\in\Z}\psi_j(x)=1\,.
  $$
 We further denote   $$\phi_j(\xi):=\sum_{k=-\infty}^j\psi_k(\xi)\,.
 $$ 
 The function $\phi_j$ is supported on  $\{\xi, ~|\xi|\le 2^{j+1}\}$.\par

 We recall the definition of the homogeneous Besov spaces $\dot{B}_{p,q}^s(\R^n)$  and  homogeneous Triebel-Lizorkin spaces $\dot{F}_{pq}^s(\R^n)$ in terms of the above dyadic decomposition.
  \begin{Definition}
  Let $s\in\R$,  $0< p,q\le\infty\,.$ For $f\in{\cal{S}}^\prime (\R^n)$ we set 
   \begin{equation}
   \left.\begin{array}{ll}
 \ds \|u\|_{\dot{B}_{p,q}^s(\R^n)}=\left(\sum_{j=-\infty}^{\infty}2^{jsq}\|{\cal{F}}^{-1}[\psi_j {\cal{F}}[u]]\|_{L^{p}(\R^n)}^q\right)^{1/q}&~~\mbox{if $q<\infty$}\\[5mm]
  \ds  \|u\|_{\dot{B}_{p,q}^s(\R^n)}=\sup_{j\in \Z} 2^{js}\|{\cal{F}}^{-1}[\psi_j {\cal{F}}[u]]\|_{L^{p}(\R^n)}&~~\mbox{if $q=\infty$}\end{array}\right.
    \end{equation}
    When $p,q<\infty$ we also set
   $$ 
   \|u\|_{\dot{F}_{p,q}^s(\R^n)}=\lf\|\left(\sum_{j=-\infty}^{\infty}2^{jsq}|{\cal{F}}^{-1}[\psi_j {\cal{F}}[u]]|^q\right)^{1/q}\rg\|_{L^p}\,.
   $$
   \hfill $\Box$
    \end{Definition}
    The space of all tempered distributions $u$ for which the quantity $\|u\|_{\dot{B}_{p,q}^s(\R^n)}$ is finite is called the homogeneous Besov space with indices 
    $s,p,q$ and it is denoted by $\dot{B}_{p,q}^s(\R^n)$. The space of all tempered distributions $f$ for which the quantity $\|f\|_{\dot F_{p,q}^s(\R^n)}$ is finite is called the homogeneous
    Triebel-Lizorkin space with indices 
    $s,p,q$ and it is denoted by $\dot F_{p,q}^s({\R}^n)\,.$ 
A classical result says \footnote{See for instance \cite{Gra1}} that $\dot{W}^{s,q}(\R^n)=\dot{B}^s_{q,2}(\R^n)=\dot{F}_{q,2}^s(\R^n)$\,.  
    
  Finally we denote    ${\cal{H}}^1(\R^n)$  the homogeneous Hardy Space in $\R^n$. A less classical results \footnote{ See for instance \cite{Gra2}.} asserts that  ${\cal{H}}^1(\R^n)\simeq \dot F^0_{2,1}$ thus we have  
  $$
   \|u\|_{{\cal{H}}^{1}(\R^n)}\simeq \int_{\R}\left(\sum_j |{\cal{F}}^{-1}[\psi_j {\cal{F}}[u]]|^2\right)^{1/2}dx\,.
 $$
  \par
  
  We recall that in dimension $n=1$,  the space  $\dot{H}^{1/2}(\R)$ is continuously embedded in the Besov space $ \dot{B}^0_{\infty,\infty}(\R)$.
 More precisely we have 
 \begin{equation} 
 \dot{H}^{1/2}(\R)\hookrightarrow BMO({\R})\hookrightarrow \dot{B}^0_{\infty,\infty}(\R)\,,\end{equation}
 where $BMO({\R})$ is the space of bounded mean oscillation dual to ${\cal{H}}^1(\R^n)$   (see for instance  \cite{RS}, page 31).
 \par 
 The $s$-fractional  Laplacian of a function  $u\colon\R^n\to\R$ is defined as a pseudo differential operator of symbol $|\xi|^{2s}$ :
 \begin{equation}\label{fract}
 \widehat {\Delta^{s}u}(\xi)=|\xi|^{2s} \hat u(\xi)\,.
 \end{equation}
 In the case where $s=1/2 $, we can write ${\Delta}^{1/2}u= -{\cal{R}}(\nabla u)$ where ${\cal{R}} $ is  Fourier multiplier of symbol
 $\displaystyle{\frac{i}{|\xi|}\sum_{k=1}^n\xi_k}:$
 $$
 \widehat {{\cal{R}}X}(\xi)=\frac{1}{|\xi|}\sum_{k=1}^n i\xi_k\hat {X_k}(\xi)$$
 for every $X\colon \R^n\to\R^n$\,, namely $ {\cal{R}}={\Delta}^{-1/2}{\rm div}\,.$\par
   
  We denote by $B_r(\bar x)$ the ball of radius $r$ and centered at $\bar x$. If $\bar x=0$ we simply write
  $B_r$\,. If $x,y\in\R^n,$ $x\cdot y$ denote the scalar product between $x,y$\,.
   
  For every function $u\colon\R^n\to\R$ we denote by $M(u)$ the maximal function of $u$, namely
 \begin{equation}\label{maxf}
 M(u)=\sup_{r>0,\ x\in\R^n}|B(x,r)|^{-1}\int_{B(x,r)}|u(y)|dy\,.
 \end{equation}

\section{Regularity of nonlocal Schr\"odinger type systems}\label{regularitySchr}
In this Section we prove Theorem \ref{regschr}. The proof is based in particular on the {\em localization estimates} established in Appendix \ref{LocEnergy}
as well on the {\em 3 commutator estimates} \rec{zz7tris} and \rec{zz9bis}\,.\par

{\bf Proof of theorem~\ref{regschr}.}
 
Let $\rho>0$ be such that $\|\11_{B(0,\rho)}\Omega\|_{L^2}\le \varepsilon_0$, with
$\varepsilon_0$ small enough. We decompose $\Om$ as follows
$\Omega_1=\11_{B(0,\rho)}\Omega$ and $\Omega_2=(1-\11_{B(0,\rho)})\Omega$\,.
\par
 Let $P\in{ \dot {H}}^{1/2}(\R, so(m))$  given by Theorem \ref{P}. We have
\begin{equation}\label{eqpv}
\Delta^{1/4} (Pv)=\left[P\Omega_1 P^{-1}-\Delta^{1/4}P P^{-1}\right]Pv+N(P,v)\end{equation}
where $N$ is the operator defined in lemma~\ref{opTintr}.\par

Since $P$ satisfies \rec{cond}(i)  we have
\begin{eqnarray}\label{sim}
P\Omega P^{-1}-\Delta^{1/4}P P^{-1}&=&-\frac{\Delta^{1/4} P P^{-1}+P\Delta^{1/4} P^{-1}}{2}\\
&=&- symm(\Delta^{1/4} P P^{-1})\nonumber\,.
\end{eqnarray}
From Theorem \ref{commDLR3bis}  it follows that 
$$
\Delta^{1/4}[symm(\Delta^{1/4} P P^{-1})]=\Delta^{1/4}[(\Delta^{1/4} P) P^{-1}]+\Delta^{1/4}[P\,\Delta^{1/4} P^{-1}]-\Delta^{1/2}(PP^{-1})\in  {\cal{H}}^1({\R})\,.
$$ 
since $PP^{-1}=Id$.  Thus $
symm(\Delta^{1/4} P P^{-1})\in L^{2,1}(\R).$\footnote{We recall that  $v\in  {\cal{H}}^1$ implies $\Delta^{-1/4} v\in L^{2,1}$ see a footnote in the introduction.}  \par
\noindent{\bf Claim 1.} 
From Theorems \ref{commDLR} and \ref{commDLR3bis} we can deduce  the estimate
\rec{zz7tris}, which can be expressed in term of the operator $N$ as follows:
     $$
\|N(Q,v)\|_{\dot{H}^{1/2}(\R^n)}\le C\|v\|_{L^{2,\infty}(\R^n)}\ \|Q\|_{H^{1/2}(\R^n)}\,.
$$
for every $Q\in{ \dot {H}}^{1/2}(\R^n)$ and $v\in L^2(\R^n)\,.$\par

\noindent{\bf Proof of Claim 1.}
\[
\begin{array}{rl}
\|N(Q,v)\|_{\dot H^{-1/2}(\R^n)}&\ds=\sup_{\|h\|_{\dot H^{1/2}}\le 1}\int_{\R^n}N(Q,v) h dx\\[5mm]
&\ds=\sup_{\|h\|_{\dot H^{1/2}}\le 1}\int_{\R^n}v [Q(\Delta^{1/4} h)-\Delta^{1/4}(Qh)+(\Delta^{1/4}Q)h] dx\\[5mm]
&\ds=\sup_{\|h\|_{\dot H^{1/2}}\le 1}\int_{\R^n} v \Delta^{-1/4}(R(Q,h))dx
\end{array}
\]
And using  Theorem \ref{commDLR3bis} we obtain
\[
\begin{array}{rl}
\|N(Q,v)\|_{\dot H^{-1/2}(\R^n)}&\lesssim \|h\|_{\dot H^{1/2}}\|v\|_{L^{2,\infty}}\|Q\|_{\dot H^{1/2}}\\[5mm]
&\lesssim \|v\|_{L^{2,\infty}}\|Q\|_{\dot H^{1/2}}\,.
\end{array}
\]
which concludes the proof of claim 1. \hfill $\Box$

\medskip

We set now $w=Pv$ and  $\omega= -symm(\Delta^{1/4} P P^{-1})$ and rewrite equation \rec{eqpv} as follows
\begin{equation}\label{eqpw}
\Delta^{1/4} w=\omega\, w+N(P,P^{-1} w)+\Omega_2 P^{-1}w\,.
\end{equation}
where by construction $\|\omega\|_{L^{2,1}}, \|P\|_{\dot H^{1/2}}\le \varepsilon_0$\,.

\medskip

\noindent{\bf Claim 2 :} {\it There exists $q>2$ such that $v\in L^q_{loc}({\R})$.}

\medskip

In order to establish the claim 2, a ''natural approach'' would be to try to prove directly a {\it Morrey decrease} for the $L^2$ norm of $v$ (or equivalently
the $L^2$ norm of $w$), that is an estimate of the form
\[
\sup_{x_0\in B(0,\rho/4)\, ,\, r>0}r^{-\beta}\left[\int_{B(x_0,r)}|v|^2\ dx\rg]^{1/2}<+\infty\,.
\]
for some positive constant $\beta>0$. We however failed to work directly with the $L^2$ norm. We have been instead  more successful in working with it's weak version : the $L^{2,\infty}-$norm.
Precisely we are going to establish the following bound 
\[
\sup_{x_0\in B(0,\rho/4)\, ,\, r>0}r^{-\beta}\ \|w\|_{L^{2,\infty}(B(x_0,r))}<+\infty\,.
\]
Let $x_0\in B(0,\rho/4)$ and $r\in(0,\rho/8)$. We argue by duality and multiply (\ref{eqpw}) by $\phi$ which is given as follows. 
Let $g\in L^{2,1}(\R)$, with $\|g\|_{L^{2,1}}\le 1$ and set $g_{r\alpha}=\11_{B(x_0,r\alpha)}g$, with $0<\alpha<1/4$ and  $\phi=\Delta^{-1/4}(g_{r\alpha})\in L^{\infty}(\R)\cap \dot H^{1/2}(\R)$\,.
We take the scalar product of both sides of equation \rec{eqpw} with $\phi$ and we integrate.

\medskip

\noindent{\bf Left hand side of the equation \rec{eqpw}:}
\begin{equation}
\label{zz10}
\begin{array}{rcl}
\ds\sup_{\|g\|_{L^{2,1}}\le 1} \int_{\R} \phi\ \Delta^{1/4} w dx&=&\ds\sup_{\|g\|_{L^{2,1}}\le 1} \int_{\R} g_{r\alpha} w dx\\[5mm]
&=&\ds\|w\|_{L^{2,\infty}(B(x_0,r\alpha))}\,.
\end{array}
\end{equation}
{\bf Right hand side of the equation \rec{eqpw}:}\par
We apply Lemmae \ref{loc1} , \ref{loc3} and \ref{loc4} and we respectively obtain in one hand
\begin{equation}
\label{zz11}
\begin{array}{rcl}
 \ds\int_{\R} \phi\ \omega\, w dx&\le&\ds \|\omega\|_{L^{2,1}}\ \|g\|_{L^{2,1}}\ \|w\|_{L^{2,\infty}(B(x_0,r))}\\[5mm]
&&\ds+\sum_{h=-1}^{+\infty} 2^{-h/2}\alpha^{1/2}\|\omega\|_{L^{2,1}}\ \|g\|_{L^{2,1}}\ \|w\|_{L^{2,\infty}(B(x_0,2^{h+1}r)\setminus B(x_0,2^{h-1}r)}
\\[5mm]
&\lesssim& \ds\varepsilon_0  \|w\|_{L^{2,\infty}(B(x_0,r))}+
\alpha^{1/2}\sum_{h=-1}^{+\infty} 2^{-h/2}\|w\|_{L^{2,\infty}(B(x_0,2^{h+1}r)\setminus B(x_0,2^{h-1}r))}\,.
\end{array}
\end{equation}
In the other hand
 \begin{equation}
 \label{zz12}
 \begin{array}{rcl}
 \ds\int_{\R} \phi\ N(P,P^{-1} w) dx&\le&\ds \varepsilon_0  \|w\|_{L^{2,\infty}(B(x_0,r)}\\
 &+&\ds C\alpha ^{1/2}\ \sum_{h=1}^{+\infty} 2^{-h/2}\|w\|_{L^{2,\infty}(B(x_0,2^{h+1}r)\setminus B(x_0,2^{h-1}r))}\,,
 \end{array}
 \end{equation}
 and finally
 \begin{equation}\label{omegalocbis}
 \int_{\R}\Omega_2 P^{-1}\,w\ \phi  dx \le C {\alpha}^{1/2}r^{1/2} \,.
 \end{equation}
 Thus combining (\ref{zz10})...(\ref{zz12}) we get
\begin{eqnarray}\label{indform}
\|w\|_{L^{2,\infty}(B(x_0,r\alpha))}&\lesssim&   \varepsilon_0  \|w\|_{L^{2,\infty}(B(x_0,r)}\\
&+& \alpha ^{1/2}\sum_{h=1}^{+\infty} 2^{-h/2}\ \|w\|_{L^{2,\infty}(B_{2^{h+1}r(x_0)}\setminus B_{2^{h-1}r}(x_0))}+  {\alpha}^{1/2}r^{1/2}\,.\nonumber
\end{eqnarray}

If $\alpha$ and $\varepsilon$ are small enough the formula \rec{indform} implies that for all
$x_0\in B(0,\rho/4)$ and $0<r<\rho/8$ we have
 $
 \|w\|_{L^{2,\infty}(B(x_0,r))}\le C r^{\beta}\,,$  for some  $\beta\in(0,1/2)\,.$ Since $P\in L^{\infty}$, this implies that
 \begin{equation}
 \label{lqcond}
\sup_{ {{\displaystyle{\mathop{\scriptstyle{x_0\in B(0,\rho/4)}}_{r>0}}}}} r^{-\beta}\ \int_{B(x_0,r)}|\Delta^{1/4} v| dx<+\infty\,.
\end{equation}
 Theorem 5.1 in \cite{AD} yields that $v\in L^q_{loc}(\R)$ for some $q>2$ which finishes the proof of claim 2.
 
 \medskip

\noindent{\bf Claim 3:} {\it $v\in L_{loc}^p(\R)$ for every $p>2$}.

\medskip

 To this end we consider again 
  $\rho>0$  such that $\|\11_{B(0,\rho)}\Omega\|_{L^2}\le \varepsilon_0$, with
$\varepsilon_0$ small enough. We decompose $\Om$ as follows. Let 
$\Omega_1=\11_{B(0,\rho)}\Omega$ and $\Omega_2=(1-\11_{B(0,\rho)})\Omega$. We consider an arbitrary $q>2$ such that $v\in L^q_{loc}$.

\medskip

Let  $x_0\in B(0,\rho/4)$, $r\in(0,\rho/8)$,
 $g\in L^{\frac{q}{q-1}}(\R)$, with $\|g\|_{L^{\frac{q}{q-1}}}\le 1$ and set $g_{r\alpha}=\11_{B(x_0,r\alpha)}g$, with $0<\alpha<1/4$ and  $\phi=\Delta^{-1/4}(g_{r\alpha})$\,.
We write the equation \rec{eqsch} as follows
\begin{eqnarray}\label{provv}
\ds\Delta^{1/4} v&=&\ds\Omega_1\11_{B(x_0,r/2)} v+\sum_{h=0}^{+\infty} \Omega_1\11_{B(x_0,2^{h+1}r)\setminus B(x_0,2^{h-1}r)} v\nonumber\\
&+&\Omega_2 v\,.
\end{eqnarray}
We take the scalar product of the equation \rec{provv} with $\Delta^{-1/4}(g_{r\alpha})$ and integrate\,.
By arguing as above, one gets
 \begin{equation}\label{lqcond3}
\sup_{ {{\displaystyle{\mathop{\scriptstyle{x_0\in B(0,\rho/4)}}_{r>0}}}}}r^{-\gamma} \lf[\int_{B(x_0,r)}| v|^q dx\rg]^{1/q}<+\infty\end{equation}
with $0<\gamma<1/4$ independent on $q\,.$
Thus by injecting \rec{lqcond3} in the equation \rec{eqsch} we obtain for the same $\gamma>0$ independent of $q$  
 \begin{equation}\label{lqcond4}
\sup_{ {{\displaystyle{\mathop{\scriptstyle{x_0\in B(0,\rho/4)}}_{r>0}}}}} r^{-\gamma}\ \|\Delta^{1/4}v\|_{L^{2q/(q+2)}B(x_0,r)}\ dx<+\infty\,.\end{equation}
Theorem 3.1 in \cite{AD} yields that $v\in L^{\tilde q}_{loc}$, with $\tilde q>q\,.$ is given by
\[
\ti{q}^{-1}=q^{-1}-2^{-1}[\gamma^{-1}(q^{-1}+2^{-1})-1]^{-1}\,.
\]
Since $q>2$ we have
\[
 \ti{q}^{-1}<q^{-1}- \frac{2\gamma}{1-4\gamma }\,.
\] 
 
By repeating the above arguments with $q$ replaced by $\tilde q$ one finally  gets that $v\in L^p_{loc}$ for every $p>2\,.$
This concludes the proof of theorem~\ref{regschr}.\hfill$\Box$

 \section{Construction of an optimal gauge $\mathbf{P}$ : the proof of theorem~\ref{P}.}\label{constrP}
 
 {\bf Proof of Theorem \ref{P}.}
 
 \medskip 
 
We follow the strategy of \cite{Riv2} in order to construct solutions to $Asymm(P^{-1}\,\Delta P)=\Om$ which was itself inspired
by Uhlenbeck's construction of Coulomb Gauges solving (\ref{zz4}).

\medskip

Let $2<q<+\infty$ and denote $1<q'<2$ to be the conjugate of $q$ : $(q)^{-1}+(q')^{-1}=1$. We consider
$$
{\cal{U}}_{\varepsilon}^q=\lf\{\Omega\in L^q(\R,so(m))\cap L^{q'}(\R,so(m)):~\int_{\R}|\Omega|^2 dx\le\varepsilon\rg\}\,.
$$

\noindent {\bf Claim:} {\it There exist $\varepsilon>0$ small enough and $C>0$ large enough such that
\[
{\cal{V}}_{\varepsilon,C}^q:=\lf\{
\begin{array}{c}
\ds \Omega\in {\cal{U}}_{\varepsilon}^q~:  \mbox{there exits $P$ satisfying \rec{cond}\ (i)-(ii)}\\[5mm]
\ds\mbox{  and }\quad\int_{\R}|\Delta^{1/4} P|^q dx\le C\int_{\R}|\Omega|^q dx
\end{array}
\rg\}
\]
is open and closed in ${\cal{U}}_{\varepsilon}^q$ and thus
${\cal{V}}_{\varepsilon}^q\equiv {\cal{U}}_{\varepsilon}^q$ (${\cal{U}}_{\varepsilon}^q$ being path  connected)}\,.
\par
\noindent {\bf Proof of the claim }\par
We first observe that ${\cal{V}}_{\varepsilon,C}^q\ne \emptyset$, ($0\in {\cal{V}}_{\varepsilon,C}^q$)\,.

\noindent{\bf Step 1:} {\it  For any $\ep>0$  and $C>0$, ${\cal{V}}_{\varepsilon,C}^q$ is closed in $L^q\cap L^{q'}(\R,so(m)).$}

\medskip

   Let $\Omega_n\in {\cal{V}}_{\varepsilon,C}^q$ such that
$\Omega_n\to \Omega_{\infty}$ in the norm $L^{q}\cap L^{q'}$, as $n\to +\infty$ and let $P_n$ be a solution
of
 \[
\begin{array}{l}
\ds P^{-1}_n\Delta^{1/4} P_n  -\Delta^{1/4} P^{-1}_nP_n=2\Omega_n\\[5mm]
\ds \int_{\R}|\Delta^{1/4}P_n|^2 dx\le C_0\int_{\R}|\Omega_n|^2 dx\,,
\end{array}
\]
Since $\Omega_n\to \Omega_{\infty}$ in the norm $L^{q}\cap L^{q'}$ and $
\int_{\R}|\Omega_n|^2 dx\le\varepsilon$, we can pass to the limit in this inequality and we have
\be
\label{zz13}
\int_{\R}|\Omega_{\infty}|^2 dx\le \varepsilon
\ee
which implies that $\Om_\infty\in {\cal{U}}_{\ep}$.

One can extract a subsequence $P_{n^{\prime}}\rightharpoonup P_{\infty}$ in $\dot H^{1/2}$. By Rellich-Kondrachov Theorem we also have 
$P_{n^{\prime}}\to P_{\infty}$ in $L^2_{loc}$ and hence $P_{\infty}\in SO(m)$ a.e. Thus $P_\infty\in \dot{H}^{1/2}({\R},SO(m))$ and the lower semi-continuity of the $\dot H^{1/2}$ and $\dot W^{1/2,q}$ norms implies that
\be
\label{zz14}
\begin{array}{cc}
 &\ds\int_{\R}|\Delta^{1/4}P_{\infty}|^2 dx\le C_0\int_{\R}|\Omega_{\infty}|^2 dx\\[5mm]
\mbox{ and } &\ds\int_{\R}|\Delta^{1/4}P_{\infty}|^q dx\le C_0\int_{\R}|\Omega_{\infty}|^q dx\,.
\end{array}
\ee
We have
\[
P^{-1}_n\Delta^{1/4} P_n  -\Delta^{1/4} P^{-1}_nP_n\to P^{-1}_{\infty}\Delta^{1/4} P_{\infty}  - \Delta^{1/4} P^{-1}_{\infty}P_{\infty}\quad\quad\mbox{ in }{\cal{D}}^{\prime}({\R})\,.
\]
Since
$P^{-1}_n\Delta^{1/4} P_n  -\Delta^{1/4} P^{-1}_nP_n=\Omega_n\to \Omega_{\infty}$ in ${\cal{D}}^{\prime}$ as well, we deduce  that
\be
\label{zz15}
P^{-1}_{\infty}\Delta^{1/4} P_{\infty}  - \Delta^{1/4} P^{-1}_{\infty}P_{\infty}=\Omega_{\infty}\quad\quad\mbox{ a.e. }
\ee
and combining (\ref{zz13}), (\ref{zz14}) and (\ref{zz15}) we deduce that $\Om_\infty\in{\cal{V}}_{\varepsilon,C}^q$ which  concludes the proof of Step 1.

\medskip

\noindent{\bf Step 2:} {\it For $\ep>0$ small enough and $C>0$ large enough ${\cal{V}}_{\varepsilon,C}^q$ is open.}

\medskip

For every $P_0\in \dot W^{1/2,q}(\R, SO(m))\cap{ \dot {H}}^{1/2}({\R},SO(m))$ we introduce the map 
 \[
 \begin{array}{l}
\ds F^{P_0}\ :\ \dot W^{1/2,q}\cap \dot{W}^{1/2,q'}(\R,so(m))\longrightarrow  L^q\cap L^{q'}(\R,so(m))\\[5mm]
\ds\quad\quad\quad U \longrightarrow  (P_0\exp U)^{-1}\Delta^{1/4}(P_0\exp U)-\Delta^{1/4} (P_0\exp U)^{-1} (P_0\exp U)\,.
\end{array}
\]
We claim first  that $F^{P_0}$ is a $C^1$ map between  the two Banach spaces $W^{1/2,q}\cap \dot{W}^{1/2,q'}(\R,so(m))$ and $L^q\cap L^{q'}(\R,so(m))$
\begin{itemize}
\item[i)] Since $\dot{W}^{1/2,q}$ for $q>2$ embedds continuously in $C^0$, the map $V\rightarrow exp\,(V)$ is clearly smooth from $\dot{W}^{1/2,q}\cap \dot{W}^{1/2,q'}({\R},so(m))$ into 
$\dot{W}^{1/2,q}\cap \dot{W}^{1/2,q'}({\R},SO(m))$.
\item[ii)] The operator $\Delta^{1/4}$ is a smooth linear map from $\dot{W}^{1/2,q}\cap \dot{W}^{1/2,q'}({\R},M_m({\R}))$ into
$L^q\cap L^{q'}({\R}, M_m({\R}))$.
\item[iii)] Since again $\dot{W}^{1/2,q}$ embedds continuously in $L^\infty$ - $\dot{W}^{1/2,q}\cap \dot{W}^{1/2,q'}$ is an algebra -
the following map
\[
\begin{array}{rcl}
\ds \Pi\ :\ \dot{W}^{1/2,q}\cap\dot{W}^{1/2,q'}({\R},M_n({\R}))\times L^q\cap L^{q'}({\R},M_n({\R}))&\longrightarrow & L^q\cap L^{q'}({\R},M_n({\R}))\\[5mm]
(A,B)\quad&\longrightarrow&\quad A\,B
\end{array}
\]
is also smooth.

\end{itemize}
Now we show that $d F_0^{P_0}=L^{P_0}$ \footnote{In order to define $L^{P_0}$ as a map from $\dot{W}^{1/2,q}\cap \dot{W}^{1/2,q'}$ into $L^q\cap L^{q'}$ we recall again that 
we make use of the embedding $\dot W^{1/2,q}(\R)\hookrightarrow L^\infty(\R)$ if $q>2$ (see for instance \cite{RS}, pag 33).} 
\[
\begin{array}{l}
\ds L^{P_0}(\eta) : = -\eta\ P_0^{-1}\Delta^{1/4} P_0+\Delta^{1/4}(\eta\ P_0^{-1})P_0\\[5mm]
\ds\quad\quad\quad\quad\quad+ P_0^{-1}\Delta^{1/4}(P_0\eta)-\Delta^{1/4}P_0^{-1} P_0\eta\,.
\end{array}
\]
\noindent$\bullet$ {\bf Differentiability of $F^{P_0}$ at $U=0$~:}\par
\[
\begin{array}{l}
\ds\lf\|F^{P_0}(\eta)-F^{P_0}(0)- L^{P_0}\cdot\eta\rg\|_{L^q\cap L^{q'}}
=\lf\|F^{P_0}(\eta)-F^{P_0}(0)+\eta P_0^{-1}\Delta^{1/4} P_0\rg.\\[5mm]
\ds\quad\quad\quad\lf.-\Delta^{1/4}(\eta P_0^{-1})P_0- P_0^{-1}\Delta^{1/4}(P_0\eta)+\Delta^{1/4}P_0^{-1} P_0\eta\rg\|_{L^q\cap L^{q'}}
\end{array}
\]
First of all we estimate
\be
\label{zz16}
\begin{array}{l}
\ds\lf\| (P_0\exp (\eta))^{-1}\Delta^{1/4}(P_0\exp U\eta)-P_0^{-1}\Delta^{1/4} P_0+\eta P_0^{-1}\Delta^{1/4} P_0-P_0^{-1}\Delta^{1/4}(\eta P_0) \rg\|_{L^q\cap L^{q'}}\\[5mm]
\ds\quad \le \lf\|\Delta^{1/4}(P_0)\rg\|_{L^q\cap L^2}\lf\|(P_0\exp(\eta))^{-1}-P_0^{-1}+\eta (P_0)^{-1}\rg\|_{L^\infty}\\[5mm]
\ds\quad\quad+\lf\|(P_0\exp(\eta))^{-1}\rg\|_{L^\infty}\ \lf\|\Delta^{1/4}(P_0\exp(\eta))-\Delta^{1/4}(P_0)-
\Delta^{1/4}(P_0 \eta)\rg\|_{L^q\cap L^{q'}}\\[5mm]
\quad\quad\ds+\lf\|\Delta^{1/4}(P_0 \eta)\rg\|_{L^q\cap L^{q'}}\ \lf\|P_0\exp(\eta)-P_0\rg\|_{L^\infty}\\[5mm]
\quad\le C\ o(\|\eta\|_{\dot{W}^{1/2,q}(\R)})
\end{array}
\ee
 The estimate of
 $$
 \lf\| (P_0{\exp  \eta})^{-1}\Delta^{1/4}(P_0\exp{(\eta)})-P_0^{-1}\Delta^{1/4} (P_0)
-P_0^{-1}\Delta^{1/4}(P_0\eta)+\Delta^{1/4}P_0^{-1} P_0\eta\rg\|_{L^q\cap L^{q'}}\,.
$$
is analogous. Hence we have proved that $dF_0^{P_0}:=L^{P_0}$.\par\bigskip

\noindent$\bullet$ {\bf $d_{0} F^{P_0}$ is an isomorphism from $\dot{W}^{1/2,q}\cap \dot{W}^{1/2,q'}(\R,so(m))$ into $L^q\cap L^{q'}(\R,so(m))$ :}\par

\medskip

Precisely we prove the following lemma.
\begin{Lemma}\label{diff}
There exists $\varepsilon>0$ such that if  $\Omega_0\in {\cal{V}}_{\varepsilon,C}^q$   and
$P_0$ is solution of \rec{cond}(i), then  for every $\omega\in L^q\cap L^{q'}(\R,so(m))$ there exists
a unique $\eta\in \dot W^{1/2,q}\cap \dot{W}^{1/2,q'}(\R,so(m))$ such that
\begin{eqnarray}\label{etaomega}
&&\omega= -\eta P_0^{-1}\Delta^{1/4} P_0+\Delta^{1/4}(\eta P_0^{-1})P_0+ P_0^{-1}\Delta^{1/4}(P_0\eta)-\Delta^{1/4}P_0^{-1} P_0\eta\,
\end{eqnarray}
and
$$
\|\eta\|_{\dot{W}^{1/2,q}\cap \dot{W}^{1/2,q'}}\le C\ \|\omega\|_{L^q\cap L^{q'}}\,.
$$
\hfill$\Box$
\end{Lemma}
{\bf Proof of Lemma \ref{diff}.}
We first observe that since $\Omega_0\in {\cal{V}}_{\varepsilon,C}^q$, then
\begin{eqnarray}
\label{zz17}
\int_{\R}|\Delta^{1/4} P_0|^2 dx&\le & C\int_{\R}|\Omega_0|^2 dx\le C\ \varepsilon\\[5mm]
\mbox{and }\quad\int_{\R}|\Delta^{1/4} P_0|^q dx&\le& C\int_{\R}|\Omega_0|^q dx\,.
\end{eqnarray}
{\bf Claim 1.}  Let $1<r<2$.  {\it $L^{P_0}$ is an isomorphism between 
$\dot W^{1/2,r}(\R,so(m))$ and $L^r$}, namely 
      for any $\omega\in L^r(\R,so(m))$ there exists a unique $\eta\in \dot W^{1/2,r}(\R,so(m))$ solution 
      to $L^{P_0}(\eta)=\om$ and
$$
\|\eta\|_{\dot W^{1/2,r}}\le C\ \|\omega\|_{L^r}\,
$$
for $C>0$\,.\par
We rewrite the equation 
\rec{etaomega} in the following way
\begin{eqnarray}\label{etaomegabis}
\omega= 2\Delta^{1/4} \eta -2\eta P_0^{-1}\Delta^{1/4}P_0-2\Delta^{1/4}P_0^{-1} P_0\eta\\[5mm]
\quad\quad+Q(\eta,P_0)-Q^{t} (\eta,P_0)\,,\nonumber
\end{eqnarray}
where 
\begin{eqnarray}\label{Q}
&&
Q(\eta,P_0)=\Delta^{1/4}(\eta P_0^{-1})P_0+\eta\ P_0^{-1}\Delta^{1/4} P_0-\Delta^{1/4}\eta\,.
\end{eqnarray}
From Lemma \ref{lemmaprel1bis} and Lemma \ref{lemmaprel2bis} it follows that 
\begin{equation}\label{Qeta}
\|Q(\eta,P_0)\|_{L^r}\le C\,\|\eta\|_{\dot W^{1/2,r}}\ \left(\|P_0\|_{\dot H^{1/2}}+\|P_0\|^2_{\dot H^{1/2}}\|P_0\|_{L^{\infty}}\right)\,.
\end{equation}
Since $2^{-1}+(2-r)\,(2r)^{-1}=r^{-1}$, by applying H\"older Inequality we get
\begin{equation}
\|\eta\ P_0^{-1}\Delta^{1/4} P_0\|_{L^r}\le \|\eta\|_{L^{{2r}/(2-r)}}\|P_0^{-1}\Delta^{1/4}P_0\|_{L^2}\,.
\end{equation}
Thus, since   $\dot W^{1/2,r}(\R,so(m))\hookrightarrow L^{\frac{2r}{2-r}}\,,$ we also have
\begin{equation}\label{P0eta}
\|\eta P_0^{-1}\Delta^{1/4} P_0\|_{L^r}\le C\, \|\eta\|_{\dot W^{1/2,r}}\|P_0^{-1}\Delta^{1/4}P_0\|_{L^2}\,.
\end{equation}
We consider the following map
$
H^{P_0} \colon \dot W^{1/2,r}(\R,so(m))\to L^r(\R,so(m))$,
\begin{eqnarray*}
H^{P_0}(\eta)&=&-2\,\eta P_0^{-1}\Delta^{1/4}P_0-2\Delta^{1/4}P_0^{-1} P_0\,\eta+ Q(\eta, P_0)-Q^t(\eta,P_0)\,.
\end{eqnarray*}
From \rec{Qeta}  and \rec{P0eta}, it follows that there exists a constant $C>0$ (independent of $P_0$) such that
$$
\|H^{P_0}(\eta)\|_{L^r}\le C\,\|\eta\|_{\dot W^{1/2,r}}\left[\|P_0\|_{\dot H^{1/2}}+\|P_0\|_{\dot H^{1/2}}^2\|P_0\|_{L^{\infty}}\right]\,.
$$
Because of (\ref{zz17}),
$\|P_0\|_{\dot H^{1/2}}\le (C\ \varepsilon)^{1/2}$ and hence, if $\varepsilon>0$ is small enough, 
$L^{P_0}=2\Delta^{1/4}+H_{P_0}\colon \dot W^{1/2,r}(\R,so(m))\to L^{r}(\R,so(m))$ is invertible
which proves the first claim.

\medskip

\noindent{\bf Claim 2.} {\it Let $q'<r<2$. Let $\om\in L^q\cap L^{r}$ and $\eta$ be the solution of $L^{P_0}(\eta)=\om$ then $\eta$ is in $\dot{W}^{q}\cap \dot{W}^{r}$.}

We apply Lemma \ref{lemmaprel4}  to
\begin{eqnarray*}
\Delta^{1/4}\eta-P_0^{-1}\Delta^{1/4}(P_0\eta)=\Delta^{1/4}(P_0^{-1}P_0\eta)-P_0^{-1}\Delta^{1/4}(P_0\eta)\,
\end{eqnarray*}
and we obtain
\begin{eqnarray}\label{estt}
&&
\|\Delta^{1/4}\eta-P_0^{-1}\Delta^{1/4}(P_0\eta)\|_{L^t}\le \|P_0\eta\|_{\dot W^{1/2,r}(\R,so(m))}\ \|P_0\|_ {\dot W^{1/2,q}(\R,so(m))}\\[5mm]
&&\le \|\eta\|_{\dot W^{1/2,r}}\left[\|P_0\|_{L^\infty}+\|P_0\|_{\dot H^{1/2}}\right]\ \|P_0\|_{\dot W^{1/2,q}(\R,so(m))}\,,\nonumber
\end{eqnarray} 
where $t$ is given by  $\frac{1}{t}=\frac{1}{q}+\frac{2-r}{2r}$.
In a similar way we   have
\begin{eqnarray}\label{estt2}
&&\|\Delta^{1/4}\eta- \Delta^{1/4}( \eta\, P^{-1}_{0})P_0\|_{L^t}\le \|\eta\|_{\dot W^{1/2,r}}\left[\|P_0\|_{L^\infty} +\|P_0\|_{\dot H^{1/2}}\right] \ \|P_0\|_{\dot W^{1/2,q}(\R,so(m)}\,.\nonumber
\end{eqnarray}
On the other hand we also have 
\begin{equation}
\|\eta\, P_0^{-1}\Delta^{1/4} P_0\|_{L^t} \le \|\eta\|_{L^{\frac{2r}{2-r}}}\ \|\Delta^{1/4}P_0\|_{L^q}
\end{equation}
Thus $Q(\eta,P_0), Q^t(\eta,P_0)$ and $H_{P_0}(\eta)$ are in $L^t\,.$ Thus since  $\omega\in L^q\cap L^r$, we have 
$\Delta^{1/4}\eta\in L^t$ as well. Since $q'<r<2$ and $\frac{1}{t}=\frac{1}{q}+\frac{1}{r}-\frac{1}{2}$, we have that $t>2$. $\Delta^{1/4}\eta\in L^r\cap L^t$ 
for some $r<2$ and $t>2$ implies that $\eta\in L^{\infty}$ (see for instance \cite{AD2}, pag 25)\,.\par
From the fact that $\eta\in L^\infty$ we deduce that
$\eta\,P_0^{-1}\Delta^{1/4} P_0\in L^q$ and $\Delta^{1/4}P_0^{-1} P_0\,\eta\in L^q$\,. Now we apply Lemma \ref{lemmaprel4}  respectively 
to $a=P_0\,\eta\in \dot H^{1/2}\cap L^{\infty}$,  $b=P_0^{-1}\in \dot W^{1/2,q}$ and $a=\eta P_0^{-1} $, $b=P_0$ and we  get that  $H_{P_0}(\eta)\in L^q$.
Since $\omega\in L^q\cap L^r$ we have $\Delta^{1/4}\eta\in L^q$ as well.
Moreover the following estimate holds
$$
\|\Delta^{1/4}\eta\|_{L^q}\le C\ \|\omega\|_{L^q\cap L^r}\le  C\ \|\omega\|_{L^q\cap L^{q^{\prime}}}\,,$$
which proves the claim 2.

\medskip

Combining claim 1 and claim 2 we obtain that for any $\om\in L^q\cap L^{q'}({\R},so(m))$ there exists a unique $\eta\in \dot{W}^{1/2,q}\cap \dot{W}^{1/2,q'}({\R},so(m))$ such that

\[
L^{P_0}\eta=\om\,, 
\]
and
\[
\|\eta\|_{\dot{W}^{1/2,q}\cap \dot{W}^{1/2,q'}}\le C\ \|\om\|_{L^q\cap L^{q'}}
\]
This finishes the proof of lemma~\ref{diff}. \hfill $\Box$

\noindent{\bf Proof of step 2 continued.}
We apply Implicit Function Theorem to $F^{P_0}$ and we  deduce that for every $P$  in  some  neighborhood  of
$P_0$ and $\Omega$ in  a  neighborhood  of
$\Omega_0$   (both neghborhoods having a size depending on $P_0$ and $\Om_0$ of course) the equation \rec{cond}(i) is satisfied and  
for some constant $C>0$ 
$$\|\Delta^{1/4} P\|_{L^q}\le C\|\Omega\|_{L^q}\,.$$
By possibly taking a smaller   neighborhood  of
$P_0$ we may always assume that $\int_{\R}|\Delta^{1/4} P|^2 dx\le \varepsilon<1$. 

\medskip

\noindent{\bf Step 3:}{\it  The fact that $\int_{\R}|\Delta^{1/4} P|^2 dx\le \varepsilon<1$  implies that
 $\int_{\R}|\Delta^{1/4} P|^2 dx\le C \int_{\R}|\Omega|^2 dx $\,.}
\par
 
We write
\begin{eqnarray*}
P^{-1}\Delta^{1/4} P&=&\frac{1}{2}(P^{-1}\Delta^{1/4} P-\Delta^{1/4}(P^{-1}\Delta^{1/4} P)^t)+\frac{1}{2}(P^{-1}\Delta^{1/4} P+(P^{-1}\Delta^{1/4} P)^t)\\
&=&
\frac{1}{2}(P^{-1}\Delta^{1/4} P-\Delta^{1/4} P^{-1}P)+\frac{1}{2}(P^{-1}\Delta^{1/4} P+\Delta^{1/4} P^{-1}P)\,.
\end{eqnarray*}
We set $$symm(P^{-1}\Delta^{1/4} P):=\frac{1}{2}(P^{-1}\Delta^{1/4} P+\Delta^{1/4} P^{-1}P)$$ and
$$Asymm(P^{-1}\Delta^{1/4} P):=\frac{1}{2}(P^{-1}\Delta^{1/4} P-\Delta^{1/4} P^{-1}P)\,.$$ 
 We apply Lemma \ref{lemmaprel2bis} and we get
 
\begin{eqnarray*}
&&\ds\int_{\R}|P^{-1}\Delta^{1/4} P+\Delta^{1/4} P^{-1}P|^2dx\le   C\|P^{-1}\Delta^{1/4} P\|^2_{L^2} \\ [5mm]&&\ds
\le C \|\Delta^{1/4} P\|_{L^2}\left(\|symm(P^{-1}\Delta^{1/4} P)\|_{L^2}+\|Asymm(P^{-1}\Delta^{1/4} P)\|_{L^2}\right)\,.
\end{eqnarray*}
Thus we get
$$
\|symm(P^{-1}\Delta^{1/4} P)\|_{L^2}\le 
C\varepsilon\left(\|sym(P^{-1}\Delta^{1/4} P)\|_{L^2}+\|Asymm(P^{-1}\Delta^{1/4} P)\|_{L^2}\right)\,.$$
If $\varepsilon >0$ is small enough then
$$\|symm(P^{-1}\Delta^{1/4} P)\|_{L^2}\le C\|Asymm(P^{-1}\Delta^{1/4} P)\|_{L^2}=C\|\Omega\|_{L^2}$$ which ends the
proof of Step 3.\par
\medskip
\noindent{\bf Step 4.} 
Take now  $\Omega\in L^2$ and $\int_{\R}|\Omega|^2 dx\le \varepsilon$. Let $\Omega_k\in{\cal{U}}_{\varepsilon}^q$ be such that $\Omega_k\to\Omega$ as $k\to +\infty$ in $L^2$.
By arguing as  in the proof of that ${\cal{V}}_{\varepsilon}^q$ is closed one gets that there exists
$P\in\dot H^{1/2}$ satisfying \rec{cond}(i)-(ii).~~\hfill$\Box$

%%%%%%%%%%%%%%%%%%%%%%%
 
\section{Euler Equation for Half-Harmonic Maps into Manifolds}\label{harmonic}
We consider a compact $k$ dimensional $C^2$ manifold without boundary ${\cal{N}}\subset\R^m$.   Let  $\Pi_{{\cal{N}}}$ be the 
orthogonal projection on ${\cal{N}}\,.$  
 We also consider the   {Dirichlet  energy} \rec{lagr}.
 
The weak ${1/2}$-harmonic maps are defined as critical points
of the functional \rec{lagr} with respect to perturbation of the form $\Pi_{{\cal{N}}}(u+t\phi)$, where $\phi$ is
an arbitrary compacted supported smooth map    from $\R $ into $\R^m\,.$
\begin{Definition}\label{weakhalfharm}
We say that $u\in H^{1/2}(\R,{\cal{N}})$ is a weak  ${1/2}$-harmonic map if and only if, for every
 maps $\phi\in  H^{1/2}(\R,\R^{m})\cap L^{\infty}(\R,\R^m)$     we have
\begin{equation}\label{critic}
\frac{d}{dt}{\cal{L}}(\Pi_{{\cal{N}}}(u+t\phi))_{|_{t=0}}=0\,.
\end{equation}
\hfill$\Box$
\end{Definition}
We introduce some notations.
 We denote by $\bigwedge(\R^m)$ the exterior algebra (or Grassmann Algebra) of $\R^m$ and by the symbol $\wedge$ the {\em exterior or wedge product}. 
For every $p=1,\ldots,m$,   $\bigwedge_p(\R^m)$ is the vector space of $p$-vectors  \par If $(\epsilon_i)_{i=1,\ldots,m}$ is the 
canonical orthonormal  basis of $\R^m$, then every element $v\in \bigwedge_p(\R^m)$ is written as
$v=\sum_{I}v_{I}\epsilon_{I}$ where $I=\{i_1,\ldots,i_p\}$ with $1\le i_1\le\ldots\le i_p\le m$ , $v_I:=v_{i_1,\ldots,i_p} $ and $ \epsilon_{I}=:=\epsilon_{i_1}\wedge\ldots\wedge \epsilon_{i_p}\,.$ \par
By the symbol $\res$ we denote the interior multiplication  $\res\colon \bigwedge_p(\R^m)\times \bigwedge_q(\R^m)\to\bigwedge_{q-p}(\R^m)$ defined as follows. \par
 Let   $\epsilon_I=\epsilon_{i_1}\wedge\ldots\wedge \epsilon_{i_p}$, $\epsilon_J=\epsilon_{j_1}\wedge\ldots\wedge \epsilon_{j_q},$ with $q\ge p\,.$ Then  $\epsilon_I\res \epsilon_J=0$ if $I\not\subset J$, otherwise
  $\epsilon_I\res \epsilon_J=(-1)^M\epsilon_{K}$ where $\epsilon_{K} $ is a $q-p$ vector and $M$ is the number of pairs $(i,j)\in I\times J$ with $j>i\,.$ \par 
 Finally by the symbol $\ast$ we denote the Hodge-star operator, $\ast\colon\bigwedge_p(\R^m)\to\bigwedge_{m-p}(\R^m)$,
defined by $\ast\beta=\beta\res(\epsilon_{1}\wedge\ldots\wedge \epsilon_n) $.  
For an introduction of the Grassmann Algebra we refer the reader to the first Chapter of the book by Federer\cite{fed}\,.\par
 In the sequel we denote by $P^T$ and $P^N$ respectively the tangent and the normal projection to the manifold ${\cal{N}}$.\par
They verify the following properties: $(P^T)^t=P^T,(P^N)^t=P^N$ (namely they are symmetric operators), $(P^T)^T=P^T$,   $(P^N)^N=P^N$,
$P^T+P^N=Id$, $P^NP^T=P^TP^N=0\,.$\par
We set $e=\epsilon_1\wedge\ldots\wedge\epsilon_k$ and $n=\epsilon_{k+1}\wedge\ldots\wedge\epsilon_m$\,. For avery $z\in {\cal{N}}$,
$e(z)$ and $n(z)$ give the orientation respectively of the tangent $k$-plane and the normal $m-k$-plane to $T_z{\cal{N}}\,.$\par
We observe that 
   for every $v\in\R^m$ we have
\begin{eqnarray}\label{tangpr}
P^T v&=&(-1)^{m-1}\ast ( ( v\res e) \wedge n) \,.
\end{eqnarray}
\begin{eqnarray}\label{normpr}
P^N v&=&(-1)^{k-1}\ast(e\wedge( v\res n)) \,.
\end{eqnarray}
We observe  that $P^N$and $P^T$ can be seen as matrices in $\dot{H}^{1/2}(\R ,\R^m)\cap
L^{\infty}(\R ,\R^m)$\,.\par
 
Next we    write the Euler equation associated to the functional \rec{lagr}\,.\par
\begin{Proposition} 
\label{pr-I.1}
All weak $1/2$-harmonic maps $u\in H^{1/2}(\R,{\cal{N}})$ satisfy in a weak sense
\par
i)  the equation
\begin{equation}\label{perp}
\int_{\R}({\Delta}^{1/2} u )\cdot v \,dx=0,
\end{equation}
for every $v\in H^{1/2}(\R,\R^{m})\cap L^{\infty}(\R,\R^m)$ and $v\in T_{u(x)}{\cal{N}}$ almost everywhere, or in a equivalent way\par
ii) the equation 
\begin{equation}\label{wedge}
 P^T{\Delta}^{1/2}u =0~~\mbox{in ${\cal{D}}^\prime\,,$}
 \end{equation}
 or 
\par
iii) the equation
\begin{equation}\label{euler1}
\Delta^{1/4} ( P^T \Delta^{1/4}u)=T(P^T,u)-(\Delta^{1/4} P^T)\Delta^{1/4} u\,,
\end{equation}
 \hfill $\Box$
 \end{Proposition}
    The Euler Lagrange equation (\ref{euler1}) can be  completed by the following ''structure equation'':

\begin{Proposition} 
\label{pr-I.2}
All maps in $\dot{H}^{1/2}(\R,{\cal{N}})$ satisfy the following identity
\begin{equation}\label{eqstruct}
 \Delta^{1/4} ({\cal{R}}(P^N \Delta^{1/4}u))={\cal{R}}(S(P^N,u))- (\Delta^{1/4} P^N) ({\cal{R}}\Delta^{1/4} u)\,.
\end{equation}
\hfill $\Box$
 \end{Proposition}
For the proofs of Proposition \ref{pr-I.1} and \ref{pr-I.2} we refer the reader to \cite{DLR}\,.\par\medskip
Next we see that by combining  \rec{euler1} and \rec{eqstruct} we can obtain the  new equation
 \rec{EulEq} for the vector field
$v=(P^T \Delta^{1/4}u,{\cal{R}}(P^N \Delta^{1/4}u))$ where an antisymmetric potential appears.  \par
  We introduce the following matrices
\begin{eqnarray} 
 &&
\omega_1=\frac{(\Delta^{1/4}P^T) P^T+P^T\Delta^{1/4}P^T-
\Delta^{1/4}(P^TP^T)}{2}\,,\label{omega1}\\
&&\omega_2={(\Delta^{1/4}P^T )P^N+P^T\Delta^{1/4}P^N-
\Delta^{1/4}(P^TP^N)}\,,\label{omega2}\\
&&\omega=\frac{(\Delta^{1/4} P^T) P^T- P^T\Delta^{1/4} P^T}{2}\,;\label{Omega}
\end{eqnarray}
and
  \begin{eqnarray} 
 &&
\omega_3=\frac{({\cal{R}}\Delta^{1/4}P^T )P^T+P^T\Delta^{1/4}({\cal{R}}\Delta^{1/4}P^T )-
{\cal{R}}\Delta^{1/4}(P^TP^T)}{2},\label{omega3}\\
&&\omega_4={({\cal{R}}\Delta^{1/4} P^T) P^N+P^N({\cal{R}}\Delta^{1/4} P^T) -{\cal{R}}\Delta^{1/4}(P^N P^T) } ,\label{omega4}\\
&&\omega_{{\cal{R}}}=\frac{({\cal{R}}\Delta^{1/4} P^T) P^T- P^T({\cal{R}}\Delta^{1/4} P^T)}{2}\,.\label{OmegaR}
\end{eqnarray}

We observe that Theorem \ref{commDLR}  and Theorem \ref{commDLR2}    imply respectively  that 
$\Delta^{1/4}(\omega_1)$,  $\Delta^{1/4}(\omega_2)$
and $\Delta^{1/4}(\omega_3)$, $\Delta^{1/4}(\omega_4)$ are in the homogeneous Hardy Space ${\cal{H}}^1(\R)$. Therefore $\omega_1,\omega_2,\omega_3,\omega_4\in L^{2,1}(\R)\,.$
The matrices  $\omega$ and $\omega_{{\cal{R}}}$ are  {\bf antisymmetric}.\par

\medskip

\noindent{\bf Proof of Proposition~\ref{EulEq}.} 
From Propositions \ref{pr-I.1} and \ref{pr-I.2} it follows that $u$ satisfies in a weak sense the equations  \rec{euler1} and \rec{eqstruct}.\par
The key point is to estimate the 
 the
  terms $(\Delta^{1/4} P^T)\Delta^{1/4} u$ and $(\Delta^{1/4} P^N) {\cal{R}}(\Delta^{1/4} u)$\par

{\bf $\bullet$  Re-writing of $(\Delta^{1/4} P^T)\Delta^{1/4} u$\,.}\par
  
\begin{eqnarray*}
(\Delta^{1/4} P^T)\Delta^{1/4} u&=&(\Delta^{1/4} P^T)(P^T\Delta^{1/4} u+P^N\Delta^{1/4} u)\\
&=&
((\Delta^{1/4} P^T) P^T)(P^T v)+((\Delta^{1/4} P^T) P^N)(P^N v)\,.
\end{eqnarray*}
Now we have
\begin{equation}\label{pt}
(\Delta^{1/4} P^T) P^T=\omega_1+\omega+\frac{\Delta^{1/4}P^T}{2}\,;
\end{equation}

 and 

 \begin{eqnarray}\label{pn}
 (\Delta^{1/4} P^T) P^N&=&( \Delta^{1/4} P^T) P^N+P^T\Delta^{1/4} P^N-\Delta^{1/4}(P^TP^N)-P^T\Delta^{1/4}P^N\nonumber\\  &=&\omega_2+P^T\Delta^{1/4} P^T\\
 &=& \omega_2+\omega_1-\omega+\frac{\Delta^{1/4}P^T}{2}\,.\nonumber
\end{eqnarray}

 Thus
 
 \begin{eqnarray}
 \frac{(\Delta^{1/4}P^T) (P^T \Delta ^{1/4} u)}{2}&=&\omega_1 ( P^T \Delta ^{1/4} u) +\omega ( P^T \Delta ^{1/4} u) \label{TT} \\ [5mm]
 \frac{(\Delta^{1/4} P^T) (P^N \Delta ^{1/4} u)}{2}&=&(\omega_1+\omega_2) (P^N \Delta ^{1/4} u) -\omega (P^N \Delta ^{1/4} u) \label{TN}\\[5mm]
 & =&{\cal{R}}(\omega_1+\omega_2){\cal{R}}(P^N \Delta ^{1/4} u)-{\cal{R}}(\omega){\cal{R}}(P^N \Delta ^{1/4} u)\nonumber\\[5mm]
 &+&F(-\omega+\omega_1+\omega_2,(P^N \Delta ^{1/4} u))\,.\nonumber
 \end{eqnarray}

{\bf $\bullet$  Re-writing of $(\Delta^{1/4} P^N) ({\cal{R}}\Delta^{1/4} u)$\,.}
\par
 We have
\begin{eqnarray*}
(\Delta^{1/4} P^N)({\cal{R}}\Delta^{1/4} u)&=&({\cal{R}}(\Delta^{1/4}P^N))(P^T(\Delta^{1/4} u)+P^N(\Delta^{1/4} u)))\\
&+&
F(({\cal{R}}(\Delta^{1/4}P^N)),\Delta^{1/4} u)\,.
\end{eqnarray*}
 
 We estimate
$({\cal{R}}\Delta^{1/4} P^N) P^T(\Delta^{1/4} u)$ and $({\cal{R}}\Delta^{1/4} P^N) P^N(\Delta^{1/4}u)$\,.
We have
\begin{eqnarray*}
({\cal{R}}\Delta^{1/4} P^N) P^T&=&-({\cal{R}}\Delta^{1/4} P^T) P^T\\[5mm]
&=&-\omega_3-\omega_{{\cal{R}}}-\frac{({\cal{R}}\Delta^{1/4} P^T)}{2}\\[5mm]
&=&
-\omega_3-\omega_{{\cal{R}}}+\frac{({\cal{R}}\Delta^{1/4} P^N)}{2}\,,
\end{eqnarray*}
and 
\begin{eqnarray*}
 ({\cal{R}}\Delta^{1/4} P^N) P^N&=&-({\cal{R}}\Delta^{1/4} P^T) P^N\pm P^T({\cal{R}}\Delta^{1/4} P^N) \\[5mm]
&=&
-[({\cal{R}}\Delta^{1/4} P^T) P^N+P^T({\cal{R}}\Delta^{1/4} P^N) -{\cal{R}}\Delta^{1/4}(P^N P^T) ]\\[5mm]
&+&
P^T({\cal{R}}\Delta^{1/4} P^N)\\[5mm]
&=&-\omega_4-\omega_3+\omega_{{\cal{R}}}+\frac{({\cal{R}}\Delta^{1/4} P^N)}{2}\,.
\end{eqnarray*}
Thus
\begin{eqnarray}
\frac{({\cal{R}}\Delta^{1/4} P^N) P^T\Delta^{1/4} u}{2}&=&-\omega_3(P^T\Delta^{1/4} u)-\omega_{{\cal{R}}}(P^T\Delta^{1/4} u)\label{NT} \\[5mm]
\frac{({\cal{R}}\Delta^{1/4} P^N) P^N\Delta^{1/4} u}{2}&=&-\omega_4(P^N\Delta^{1/4} u)-\omega_3(P^N\Delta^{1/4} u)+\omega_{{\cal{R}}}(P^N\Delta^{1/4} u)\label{NN}\\
&=&{\cal{R}}(-\omega_3-\omega_4){\cal{R}}(P^N\Delta^{1/4} u)\nonumber\\[5mm]
&+&{\cal{R}}( \omega_{{\cal{R}}}){\cal{R}}(P^N\Delta^{1/4} u)\nonumber\\[5mm]
&+&F( \omega_{{\cal{R}}}- \omega_3- \omega_4,P^N\Delta^{1/4} u)\,.\nonumber
\end{eqnarray}

By combining \rec{TT}, \rec{TN}, \rec{NT} , \rec{NN}  we obtain

\begin{eqnarray}\label{Eulerbis}
\Delta^{1/4}\left(\begin{array}{l}
P^T\Delta^{1/4} u\\ 
{\cal{R}}P^N\Delta^{1/4} u\end{array}\right)&=& \tilde \Omega_1+ \tilde \Omega_2\left(\begin{array}{l}
P^T\Delta^{1/4} u\\ 
{\cal{R}}P^N\Delta^{1/4} u\end{array}\right) \\ 
&+&
2\left(\begin{array}{cc}
- \omega &  \omega_{{\cal{R}}}\\ 
 \omega_{{\cal{R}}}  & -{\cal{R}} \omega_{{\cal{R}}} \end{array}\right)
\left(\begin{array}{l}
P^T\Delta^{1/4} u\\
{\cal{R}}P^N\Delta^{1/4} u\end{array}\right)\,,\nonumber
\end{eqnarray}
 where
 $\tilde \Omega_1 $ and  $\tilde  \Omega_2 $ are given by
$$
\tilde \Omega_1=\left(\begin{array}{c}
- 2F( -\omega+ \omega_1+ \omega_2,(P^N \Delta ^{1/4} u))+ T(P^T,u)\\
-2F({\cal{R}}(\Delta^{1/4} P^N),{\cal{R}}(\Delta^{1/4} u))-2F(\omega_{{\cal{R}}}-\omega_3-\omega_4,P^N(\Delta^{1/4} u)+{\cal{R}}(S(P^N,u))
\end{array}\right)\,.$$

$$
\tilde \Omega_2=2\left(\begin{array}{cc}
-\omega_1 & -[{\cal{R}}(\omega_1+\omega_2)+({\cal{R}}(\omega)-\omega_{{\cal{R}}})]\\
\omega_3 & -{\cal{R}}(\omega_3-\omega_4)\end{array}\right)\,.$$

The matrix
$$
\Omega=2\left(\begin{array}{cc}
- \omega &  \omega_{{\cal{R}}}\\
 \omega_{{\cal{R}}}  & - {\cal{R}} \omega_{{\cal{R}}} \end{array}\right)$$
is antisymmetric\,.
 
We observe that from the estimate \rec{estFinftybis} it follows that
   $\tilde \Omega_1\in {\cal{H}}^{-1/2}(\R,\R^{2m})$ and
 \begin{equation}\label{tildeomega1}
\|\tilde\Omega_1\|_{ H^{-1/2}(R)}\le C( \|P^N\|_{\dot H^{1/2}(\R)}+ \|P^T\|_{\dot H^{1/2}(\R)})\|\Delta^{1/4} u\|_{L^{2,\infty}}\,.
\end{equation}
  On the other hand  $\tilde  \Omega_2\in L^{2,1}(\R,{\cal{M}}_{2m})$ and
   \begin{equation}\label{tildeomega2}
\|\tilde\Omega_2\|_{ L^{2,1}(R)}\le C( \|P^N\|^2_{\dot H^{1/2}(\R)}+ \|P^T\|^2_{\dot H^{1/2}(\R)})\,.
\end{equation}
This concludes the proof of  proposition~\ref{EulEq}.\hfill$\Box$

\medskip

\noindent{\bf Proof of Theorem \ref{reghm}.} \par
From Proposition  \ref{EulEq}  it follows that $v=(P^T(\Delta^{1/4} u),{\cal{R}}(P^N(\Delta^{1/4} u)))$ solves equation 
    \rec{Eulerbis}  which  of the type \rec{eqsch} up to the term $\tilde\Omega_1$.
    Therefore  the arguments are very  similar to those  of Theorem \ref{regschr} and we give
only a sketch of proof.\par
  We aim at obtaining that  $\Delta^{1/4} u\in L^{p}_{loc}(\R)$, for all $p\ge 1\,.$    To this purpose  we take $\rho>0$ such that
$$\|\Omega\|_{L^2(B(0,\rho)}, \|P^T\|_{\dot H^{1/2}(B(0,\rho)}, \|P^N\|_{\dot H^{1/2}(B(0,\rho)}\le \varepsilon_0,$$ with $\varepsilon_0>0$ small enough. 
Let $x_0\in B(0,\rho/4)$ and $r\in(0,\rho/8)$. As in the case of  equation \rec{eqsch} we argue by duality and multiply  both sides of equation \rec{Eulerbis}   by  $\phi=\Delta^{-1/4}(g_{r\alpha})$,  with
  $g\in L^{2,1}(\R)$,  $\|g\|_{L^{2,1}}\le 1$ and  $g_{r\alpha}=\11_{B(x_0,r\alpha)}g$, with $0<\alpha<1/4$\,.    \par
It is enough to estimate the integral \begin{equation}\label{lastint}
\int_{\R}\tilde\Omega_1\left(\begin{array}{c}\Delta^{-1/4} P^T(g_{r\alpha})) \\ \Delta^{-1/4} P^N(g_{r\alpha}))\end{array}\right) dx\,.
\end{equation}
 (being the other terms already estimated in the proof of  Theorem \ref{regschr})\,.

We observe that 
\begin{equation}\label{cnorme}
\|\Delta^{1/4} u\|_{L^{2,\infty}}\lesssim\left \|\sqrt{(P^T(\Delta^{1/4} u))^2+({\cal{R}}(P^N(\Delta^{1/4} u)))^2}\right\|_{L^{2,\infty}}=\|v\|_{L^{2,\infty}}
\end{equation}
By combining Lemma \ref{loc1}, \ref{loc2}, \ref{loc5} and \ref{loc6} and the estimate \rec{cnorme} we obtain
\begin{eqnarray*}
 \rec{lastint}&\lesssim &\varepsilon_0\|\Delta^{1/4} u\|_{L^{2,\infty}}+ \alpha^{1/2}\sum_{h=1}^{+\infty} 2^{-h/2}\|\Delta^{1/4} u\|_{L^{2,\infty}(B_{2^{h+1}r}\setminus B_{2^{h-1}r})}\\
 &\lesssim &\varepsilon_0\|v\|_{L^{2,\infty}}+ \alpha^{1/2}\sum_{h=1}^{+\infty} 2^{-h/2}\|v\|_{L^{2,\infty}(B_{2^{h+1}r}\setminus B_{2^{h-1}r})}\,.
\end{eqnarray*}
Therefore $v$ satisfies an estimate of the type \rec{indform} which
  implies $
 \|v\|_{L^{2,\infty}(B(x_0,r))}\le C r^{\beta}\,,$ for  $\alpha$ and $\varepsilon_0$   small enough, for all
$x_0\in B(0,\rho/4)$ and $0<r<\rho/8$  and
   for some  $\beta\in(0,1/2)\,.$ 
 
By arguing as in Theorem \ref{regschr}   we deduce that $ v\in L^{p}_{loc}(\R)$, for all $p\ge 1\,.$  Therefore $\Delta^{1/4} u\in L^{p}_{loc}(\R)$, for all $p\ge 1$ as well.

This implies that $u\in C^{0,\alpha}_{loc}$ for some $0<\alpha<1$, since $W^{1/2,p}_{loc}(\R)\hookrightarrow C^{0,\alpha}_{loc}(\R)$ if $p>2$
(see for instance \cite{AD2}). This concludes the proof of Theorem \ref{reghm}\,.~\hfill$\Box$
 
 \appendix
 \section{Localization Estimates}\label{LocEnergy}
 The aim of this Appendix 
  is to provide 
    {localization estimates} for the terms appearing in the equation 
 \rec{eqsch} and the equation \rec{Eulerbis}\,.  
 
 \par
 For $r>0$, $h\in\Z$ and $x_0\in\R$ we set 
  $$A_{h,x_0}=B (x_0,2^{h+1})\setminus B (x_0,2^{h-1})~~\mbox{ and}~~ 
A^\prime_{h,x_0}=B (x_0,2^{h})\setminus B (x_0,2^{h-1})\,.$$ 
 % Localization 1
 
 We first localize the term $N(Q,v)$.\par
 
 \begin{Lemma}\label{loc1}
 Let $Q\in{ \dot {H}}^{1/2}(\R)\cap L^{\infty}(\R)$, $\|Q\|_{\dot H^{1/2}(\R)}\le \varepsilon_0,$ $v\in L^2(\R)$,
 $g\in L^{2,1}(\R)$,  supp $g\subset B(x_0,r\alpha)$, with $x_0\in\R$, $0<\alpha<\frac{1}{4}$, $r>0$. 
 Then we have
 \begin{eqnarray}\label{eqT}
 \int_{\R} N(Q,v)\Delta^{-1/4} g dx &\lesssim& \varepsilon_0 \|g\|_{L^{2,1}} \|v\|_{L^{2,\infty}(B(x_0,r))}\\
 &+&(\|Q\|_{\dot H^{1/2}(\R)}+\|Q\|_{L^{\infty}})\|g\|_{L^{2,1}} \alpha^{1/2}\sum_{h=1}^{+\infty} 2^{-h/2}\|v\|_{L^{2,\infty}(A_{h,x_0})}\,.\nonumber
 \end{eqnarray}
 \end{Lemma}
 {\bf Proof of Lemma \ref{loc1}.} 
  We  consider a dyadic decomposition of the unity   $\varphi_j\in C_0^\infty(\R)$    such that
  \begin{equation}\label{unity}
  \mbox{ supp$(\varphi_j) \subset B_{2^{j+1}r}(x_0)\setminus B_{2^{j-1}r}(x_0)$},~~\sum_{-\infty}^{+\infty}\varphi_j=1\,.\end{equation}
   \par
   We set
    $\chi_r:=\sum_{-\infty}^{0}\varphi_j$.   
We observe that the function $\psi=\Delta^{-1/4} g$ is in $L^{\infty}(\R)\cap{ \dot {H}}^{1/2}(\R)$.\par
We take the scalar product of $N(Q,v)$ with $\Delta^{-1/4} g$ and we integrate. We write
\begin{eqnarray*}
\int_{\R} N(Q,v) \Delta^{-1/4} g dx&=&\underbrace {\int_{\R} N(Q,\chi _r v) \Delta^{-1/4} g dx}_{(1)}\\
&+& \underbrace{\int_{\R}\sum_{h=1}^{+\infty} N(Q,\varphi_h v) \Delta^{-1/4} g dx}_{(2)}
\end{eqnarray*}
To estimate $(1)$ we use the fact that $N(Q,v)\in \dot H^{-1/2}(\R)$ and $$\|N(Q,v)\|_{\dot H^{1/2}(\R)}\lesssim
\|Q\|_{\dot H^{1/2}(\R)}\|v\|_{L^{2,\infty}(\R)}\,.$$
\begin{eqnarray*}
(1)&\le&  \|\Delta^{-1/4} g\|_{{ \dot {H}}^{1/2}(\R)}\|Q\|_{\dot H^{1/2}(\R)}\|v\|_{L^{2,\infty}}\\
&\lesssim &  \varepsilon_0 \|g\|_{L^{2,1}} \|v\|_{L^{2,\infty}(B(x_0,r)}\,.
\end{eqnarray*}
Next we spilt $(2)$ in two parts:\par
 \begin{eqnarray*}
 (2) 
 &=&\underbrace{ \sum_{k=1}^{\infty}\int_{\R} N(Q,\varphi_k v)\11_{B(x_0,r/4) }\Delta^{-1/4} g dx}_{(3)}\\
 &+&\underbrace{\sum_{k=1}^{\infty}\sum_{h=-1}^{\infty}{\int_{\R} N(Q,\varphi_k v)\11_{ A^\prime_{h,x_0}}}
  \Delta^{-1/4} g dx}_{(4)}\,.
  \end{eqnarray*}
  We observe that in $(3)$ and $(4)$ we can exchange the integral with the infinite sum (see the Appendix in \cite{DLR})\,.\par
We estimate $(3)$. We first observe that since  $\11_{B(x_0,r/4) }$ and  $\varphi_k$ have disjoint supports, we have 
$$
N(Q,\varphi_k v)\11_{B(x_0,r/4) }=\left[\Delta^{1/4}(Q\varphi_k v)-Q\Delta^{1/4}( \varphi_k v)\right])\11_{B(x_0,r/4) }\Delta^{-1/4} g dx\,.$$
Thus
\begin{eqnarray*}
&&(3)=  \sum_{k=1}^{\infty}\int_{\R} \left[\Delta^{1/4}(Q\varphi_k v)-Q\Delta^{1/4}( \varphi_k v)\right])\11_{B(x_0,r/4) }dx \\[5mm]
&&\simeq \sum_{k=1}^{\infty}\int_{\R}{\cal{F}}^{-1}(|\cdot|^{1/2})(\xi)\\
&&~~~~\left[Q(\varphi_k v)\ast( \11_{B(x_0,r/4) }\Delta^{-1/4} g)-(\varphi_k v)\ast (Q \11_{B(x_0,r/4) }\Delta^{-1/4} g)\right] d\xi\\[5mm]
&&
\lesssim 
 \sum_{k=1}^{\infty}\|\xi|^{-3/2}|\|_{L^{\infty}(B^c(0,2^hr))}
 \\[5mm]
 &&~~\left[\|Q(\varphi_k v)\ast( \11_{B(x_0,r/4) }\Delta^{-1/4} g)\|_{L^1(\R)}+\|(\varphi_k v)\ast (Q \11_{B(x_0,r/4) }\Delta^{-1/4} g)\|_{L^1}\right]\\[5mm]
&&\lesssim  
 \sum_{k=1}^{\infty} 2^{-3/2k}r^{-3/2}\left[
 2\|Q\|_{L^{\infty}}\|\varphi_k v\|_{L^1(A_{h,x_0})}\|\11_{B(x_0,r/4) }\Delta^{-1/4} g)\|_{L^1(\R)}\right]\\[5mm]
&& \lesssim  
 \sum_{k=1}^{\infty} 2^{-3/2k}r^{-3/2} 2^{k/2}r^{1/2}r^{1/2} (r\alpha)^{1/2}
\left[\|Q|_{L^{\infty}}\| v\|_{L^{\infty}(A_{h,x_0})}\|g\|_{L^{2,1}(\R)}\right]\\[5mm]
&&
\lesssim  \alpha^{1/2}  \|Q\|_{L^{\infty}(\R)}\|g\|_{L^{2,1}(\R)}
  \sum_{k=1}^{\infty}2^{-k/2}\| v\|_{L^{\infty}(A_{k,x_0})}\,.
 \end{eqnarray*}
 We estimate $(4)$\,. 
\begin{eqnarray*}
&&
(4)=\underbrace{\sum_{k=1}^{\infty}\sum_{|k-h|\le 5}{\int_{\R} N(Q,\varphi_k v)\11_{ A^\prime_{h,x_0}}}\Delta^{-1/4} g dx}_{(5)}\\
&&
=\underbrace{\sum_{k=1}^{\infty}\sum_{|k-h|\ge 5}{\int_{\R} N(Q,\varphi_k v)\11_{ A^\prime_{h,x_0}}}\Delta^{-1/4} g dx}_{(6)}\,.\end{eqnarray*}
We observe that
$$
\|\11_{ A^\prime_{h,x_0}}\Delta^{-1/4} g\|_{\dot H^{1/2}(\R)}\lesssim \|g\|_{L^{2,1}(\R)}\alpha^{1/2} 2^{-h/2} \,.
$$
Thus
\begin{eqnarray*}
(5)&\le &  \sum_{k=1}^{\infty}\sum_{|k-h|\ge 5}\| N(Q,\varphi_k v)\|_{\dot H^{-1/2}(\R)}\|\11_{ A^\prime_{h,x_0}}\Delta^{-1/4} g\|_{\dot H^{1/2}(\R)}\\
&
\lesssim&\alpha ^{1/2}\|Q\|_{\dot H^{1/2}(\R)}\|g\|_{L^{2,1}(\R)}\sum_{k=1}^{\infty} 2^{-k/2} \| v\|_{L^{2,\infty}(A_{k,x_0})}\,.
\end{eqnarray*}
In order to estimate $(6)$ we observe if $|k-h|\ge 6$ then  $\varphi_k v$ and $\11_{ A^\prime_{h,x_0}}\Delta^{-1/4} g$ have disjoint supports . Thus by arguing as in $(3)$ we get
\begin{eqnarray*}
(6)&\lesssim&\alpha^{1/2} \|Q\|_{L^{\infty}}\|g\|_{L^{2,1}(\R)} \sum_{k=1}^{+\infty} 2^{-k/2}  \| v\|_{L^{2,\infty}(A_{k,x_0})} \,.~ 
\end{eqnarray*}
This concludes the proof of Lemma \ref{loc1}\,.$\hfill\Box$
\begin{Lemma}\label{loc2}
 Let $Q\in{ \dot {H}}^{1/2}(\R)\cap L^{\infty}(\R)$,  supp $Q\subset B^c(0,\rho)$  for some $\rho>0$, $v\in L^2(\R)$, $x_0\in B(0,\rho/4)$, 
 $g\in L^{2,1}(\R)$,  supp $g\subset B(x_0,r\alpha)$, with  , $0<\alpha<1$, $0<r<\rho/8$.\par
 Then we have
 \begin{eqnarray}\label{eqT2}
 &&
 \int_{\R} N(Q,v)\Delta^{-1/4} g dx \lesssim\left( \frac{r}{\rho}\right)^{1/2}\|g\|_{L^{2,1}(\R)}\|Q\|_{{ \dot {H}}^{1/2}(\R)}\|v\|_{L^{2,\infty}(B(x_0,r))}\\
 &&+ \alpha ^{1/2}\ (\|Q\|_{L^{\infty}}+\|Q\|_{{ \dot {H}}^{1/2}(\R)})\|g\|_{L^{2,1}(\R)} \sum_{h=1}^{+\infty} 2^{-h/2}\|v\|_{L^{2,\infty}(A_{h,x_0})}\,.\nonumber \end{eqnarray}
 \end{Lemma}
 {\bf Proof of Lemma \ref{loc2}.} We write
 \begin{eqnarray*} 
 \int_{\R} N(Q,v)\Delta^{-1/4} g dx&=& \underbrace{\int_{\R} N(Q,\chi_r v)\Delta^{-1/4} g dx}_{(7)}\\
 &=& \underbrace{ \int_{\R} N(Q,(1-\chi_r) v)\Delta^{-1/4} g dx}_{(8)}\,.
 \end{eqnarray*}
  We denote by $Q_{\rho}=|B_{\rho}(0)|^{-1}\int_{B_{\rho}(0)} Q(y)dy=0$ and 
 write $Q=\sum_{h=-1}^{+\infty}\tilde \varphi_h(Q-Q_{\rho})\,,$ with supp$(\tilde\varphi_h) \subset B(0,2^{h+1}\rho)\setminus  B(0,2^{h-1}\rho)$, $\tilde\varphi$ partition of unity.\par We estimate $(7)$.\par
 \begin{eqnarray*}
 &&
 (7)=\int_{\R}N(\sum_{h=-1}^{+\infty} \varphi_h(Q-Q_{\rho}),\chi_r v)\Delta^{-1/4} g dx\\
 && =\sum_{h=-1}^{+\infty} \int_{\R} [-\varphi_h(Q-Q_{\rho})\Delta^{1/4}(\chi_r v)\Delta^{-1/4} g + \Delta^{1/4}(\varphi_h(Q-Q_{\rho})(\chi_r v)\Delta^{-1/4} g] dx\\
 &&
 =\sum_{h=-1}^{+\infty}{\cal{F}}^{-1}[|\cdot|^{1/2}](\xi)\\
 &&~~~\left[-(\chi_r v)\ast( \varphi_h(Q-Q_{\rho})\Delta^{-1/4} g)+ \varphi_h(Q-Q_{\rho})\ast( \chi_r v \Delta^{-1/4} g)\right] dx\\
 &&
 \lesssim \sum_{h=-1}^{+\infty}\||\xi|^{-3/2}\|_{L^{\infty}(B^c(0,2^h\rho))}\left[\|\chi_r v\|_{L^1}\|\varphi_h(Q-Q_{\rho})\|_{L^1}\|\Delta^{-1/4} g\|_{L^{\infty}}\right]\\
 &&
 \lesssim \|g\|_{L^{2,1}(\R)} \sum_{h=-1}^{+\infty} 2^{-h/2}\left(\frac{r}{\rho}\right)^{1/2}\|v\|_{L^{2,\infty}(B(x_0,r)}
 \|\varphi_h(Q-Q_{\rho})\|_{{ \dot {H}}^{1/2}(\R)}\\
 &&\mbox{by Lemma 4.1 in \cite{DLR}}\\
 &&
\lesssim \left( \frac{r}{\rho}\right)^{1/2}\|g\|_{L^{2,1}(\R)} \|Q\|_{{ \dot {H}}^{1/2}(\R)}\|v\|_{L^{2,\infty}(B(x_0,r))}\,.
\end{eqnarray*}
By arguing as in  $(3)$ and $(4)$ we get 
\begin{equation}
(8)\lesssim (\|Q\|_{L^{\infty}}+\|Q\|_{{ \dot {H}}^{1/2}(\R)})\|g\|_{L^{2,1}(\R)}\alpha ^{1/2}\sum_{h=1}^{+\infty} 2^{-h/2}\|v\|_{L^{2,\infty}(A_{h,x_0})}\,. 
\end{equation}
This concludes the proof of Lemma \ref{loc2}\,.~\hfill$\Box$\par
\medskip
The localization of the operator $S(Q,\Delta^{-1/4} v)$, with $v\in L^2(\R)$  is similar to that of $N(Q,v)$ and we omit it.\par\bigskip
%localization L^{2,1}
\begin{Lemma}\label{loc3}
Let $A\in L^{2,1}(\R)$, $x_0\in\R$, $r>0$,  $0<\alpha<1$ and $g\in L^{2,1}(\R)$,  supp $g\subset B(x_0,r\alpha)$. Then
\begin{eqnarray}\label{Aloc}
\int_{\R} Av \Delta^{-1/4} g  dx&\lesssim& \|A\|_{L^{2,1}}\|g\|_{L^{2,1}}\|v\|_{L^{2,\infty}(B(x_0,r))}\\
&+&\alpha^{1/2} \sum_{h=-1}^{+\infty} 2^{-h/2} \|A\|_{L^{2,1}}\|g\|_{L^{2,1}}\|v\|_{L^{2,\infty}(A_{h,x_0})}\,.\nonumber
\end{eqnarray}
\end{Lemma}
{\bf Proof of Lemma \ref{loc3}.}
 We write
 \begin{eqnarray*}
\int_{\R} Av \Delta^{-1/4} g  dx&=&\underbrace{
\int_{\R} Av\11_{B(x_0,r) }\Delta^{-1/4} g  dx}_{(9)}+\underbrace{
\sum_{h=0}^{+\infty}\int_{\R} Av\11_{A^{\prime}_{h,x_0}}\Delta^{-1/4} g  dx}_{(10)}
\end{eqnarray*}
We have
\begin{eqnarray*}
(9)&\le& \|A \Delta^{-1/4} g \|_{L^{2,1}}\|v\|_{L^{2,\infty}(B(x_0,r))}\\[5mm]
&\le& \|A\|_{L^{2,1}}\| \Delta^{-1/4} g\|_{L^{\infty}}\|v\|_{L^{2,\infty}(B(x_0,r))}\\[5mm]
&\lesssim&
 \|A\|_{L^{2,1}}\|g\|_{L^{2,1}}\|v\|_{L^{2,\infty}(B(x_0,r))}\,.
 \end{eqnarray*}
\begin{eqnarray*}
(10)&\simeq& \sum_{h=0}^{+\infty}\int_{\R}{\cal{F}}^{-1}[|\cdot|^{-1/2}](\xi) g\ast (\11_{A^{\prime}_{h,x_0}}A v)d\xi\\
&\lesssim&\sum_{h=0}^{+\infty}
\|\xi|^{-1/2}\|_{L^{\infty}(B^c(0,2^hr)}\|g\ast (\11_{A^{\prime}_{h,x_0}}A v)\|_{L^1}\\
& \lesssim& \sum_{h=0}^{+\infty} 
2^{-h/2} r^{-1/2}\|g\|_{L^1}\|\11_{A^{\prime}_{h,x_0}}A v\|_{L^1}\\
&
\lesssim&\sum_{h=0}^{+\infty} 2^{-h/2}r^{-1/2}(r\alpha)^{1/2} \|g\|_{L^{2,1}}\|A\|_{L^{2,1}}\|v\|_{L^{2,\infty}(A_{h^{\prime},x_0})}\\
&
\lesssim& \alpha^{1/2} \|g\|_{L^{2,1}}\|A\|_{L^{2,1}} \sum_{h=0}^{+\infty}2^{-h/2} \|v\|_{L^{2,\infty}(A_{h,x_0})}\,.
\, 
\end{eqnarray*}
This concludes the proof of Lemma \ref{loc3}\,.~\hfill$\Box$\par
\medskip
\begin{Lemma}\label{loc4}
Let $\Omega\in  L^2(\R,{\cal{M}}_{m\times m}(\R)$ be such that \mbox{supp $\Omega\subset B^c(0,\rho)$},
$v\in L^2(\R)$, $x_0\in B(0,\rho/4)$, 
 $g\in L^{2,1}(\R)$,  supp $g\subset B(x_0,r\alpha)$, with  , $0<\alpha<1$, $0<r<\rho/8$.\par
 Then we have
 \begin{equation}\label{omegaloc}
 \int_{\R}\Omega v\Delta^{-1/4} g dx \lesssim {(r\alpha)}^{1/2}  \|g\|_{L^{2,1}}\|\Omega\|_{L^{2}}\|v\|_{L^{2}}\,.\end{equation}
  \end{Lemma}
{\bf Proof of Lemma \ref{loc4}.}  
We use the fact that $\Omega$ and $g$ have disjoint supports.
\begin{eqnarray*}
 \int_{\R}\Omega v\Delta^{-1/4} g dx&=&\int_{\R}{\cal{F}}^{-1}(|\cdot|^{-1/2})(\xi) g\ast \Omega v d \xi\\
 &
 \lesssim& \| |x|^{-1/2}\|_{B^c(0,\rho/4)}\| g\ast \Omega v \|_{L^1}\\
 & \lesssim  &\left(\frac{4}{\rho}\right)^{1/2}\|g\|_{L^1}\|\Omega v \|_{L^1}\\
 &\lesssim &  {(r\alpha)}^{1/2}  \|g\|_{L^{2,1}}\|\Omega\|_{L^2}\|v\|_{L^{2}} \,. 
 \end{eqnarray*}
 This concludes the proof of Lemma \ref{loc4}\,.~\hfill$\Box$\par
\medskip
 Now we are going to localize the operator $F$ defined in \rec{F}\,.\par
 
  \begin{Lemma}\label{loc5}
 Let $Q\in{ L^2}(\R)\cap L^{\infty}(\R)$, $\|Q\|_{L^2(\R)}\le \varepsilon_0,$ $v\in L^2(\R)$,
 $g\in L^{2,1}(\R)$,  supp $g\subset B(x_0,r\alpha)$, with $x_0\in\R$, $0<\alpha<\frac{1}{4}$, $r>0$.\par
 Then we have
 \begin{eqnarray}\label{eqF}
 \int_{\R} F(Q,v)\Delta^{-1/4} g dx &\lesssim& \varepsilon_0 \|g\|_{L^{2,1}} \|v\|_{L^{2,\infty}(B_r(x_0)}\\
 &+&\alpha ^{1/2} (\|Q\|_{L^2(\R)}+\|Q\|_{L^{\infty}}) \|g\|_{L^{2,1}} \sum_{h=1}^{+\infty} 2^{-h/2}\|v\|_{L^{2,\infty}(A_{h,x_0})}\,.\nonumber
 \end{eqnarray}
 \end{Lemma}
 {\bf Proof of Lemma \ref{loc5}.} 
 We take the scalar product of $F(Q,v)$ with $\Delta^{-1/4} g$ and we integrate. We get
\begin{eqnarray*}
\int_{\R} F(Q,v) \Delta^{-1/4} g dx&=&\underbrace {\int_{\R} F(Q,\chi _r v) \Delta^{-1/4} g dx}_{(11)}\\
&+& \underbrace{\int_{\R}\sum_{k=1}^{+\infty} F(Q,\varphi_k v) \Delta^{-1/4} g dx}_{(12)}\,.
\end{eqnarray*}
To estimate $(11)$ we use the fact that $F(Q,v)\in \dot H^{-1/2}(\R)$ and $$\|F(Q,v)\|_{\dot H^{1/2}(\R)}\lesssim
\|Q\|_{L^2(\R)}\|v\|_{L^{2,\infty}}\,.$$
\begin{eqnarray*}
(11) &\le&  \|\Delta^{-1/4} g\|_{{ \dot {H}}^{1/2}(\R)}\|Q\|_{L^2(\R)}\|v\|_{L^{2,\infty}(B(x_0,r))}\\
&\lesssim &  \|g\|_{L^{2,1}}\|Q\|_{L^2(\R)}\|v\|_{L^{2,\infty}(B(x_0,r))}\\
&\lesssim &  \varepsilon_0\|g\|_{L^{2,1}} \|v\|_{L^{2,\infty}(B(x_0,r))}
\,.
\end{eqnarray*}
Next we spilt $(12)$ in two parts:\par
 \begin{eqnarray*}
 (12) 
 &=&\underbrace{ \sum_{k=1}^{\infty}\int_{\R} F(Q,\varphi_k v)\11_{B(x_0,r/4) }\Delta^{-1/4} g dx}_{(13)}\\
 &+&\underbrace{\sum_{k=1}^{\infty}\sum_{h=-1}^{\infty}{\int_{\R} F(Q,\varphi_k v)\11_{ A^\prime_{h,x_0}}}
  \Delta^{-1/4} g dx}_{(14)}\,.
  \end{eqnarray*}
  Estimate of $(13)$:\par
  \begin{eqnarray*}
  (13)&=&\sum_{k=1}^{+\infty}\int_{\R}F(Q,\varphi_k v)\11_{B(x_0,r/4)}\Delta^{-1/4} g dx\\[5mm]
  &=&\sum_{k=1}^{+\infty}\int_{\R}{\cal{R}}(Q){\cal{R}}(\varphi_k v)\11_{B(x_0,r/4)}\Delta^{-1/4} g dx\\[5mm]
  &\simeq& \sum_{k=1}^{+\infty}\int_{\R}{\cal{F}}^{-1}\left[\frac{\cdot}{|\cdot|}\right](\xi)(\varphi_k v)\ast(Q\11_{B(x_0,r/4)}\Delta^{-1/4} g)d\xi\\
  &\lesssim & \sum_{k=1}^{+\infty}\|\frac{1}{\xi}\|_{L^{\infty}(B^c(x_0,2^{k-1}r))}\|\varphi_k v\|_{L^1(\R)}\|Q\11_{B(x_0,r/4)}\Delta^{-1/4} g\|_{L^1(\R)}\\[5mm]
  &\lesssim& \sum_{k=1}^{+\infty} 2^{-k}r^{-1}2^{k/2}r^{1/2}r\alpha^{1/2} \|v\|_{L^{2,\infty}(A_{h,x_0})}\|Q\|_{L^{\infty}} \|g\|_{L^{2,1}}\\[5mm]
   &\lesssim&(r\alpha)^{1/2}   \|Q\|_{L^{\infty}} \|g\|_{L^{2,1}}  \sum_{k=1}^{+\infty} 2^{-k/2}\|v\|_{L^{2,\infty}(A_{h,x_0})}\,.
   \end{eqnarray*}
  
 The estimate of $(14)$ is analogous of
$(4)$ in the proof of Lemma \ref{loc2} and we omit it.~\hfill$\Box$
\par
\begin{Lemma}\label{loc6}
 Let $Q\in{ L^2}(\R)\cap L^{\infty}(\R)$,  supp $Q\subset B^c(0,\rho)$  for some $\rho>0$, $v\in L^2(\R)$, $x_0\in B(0,\rho/4)$, 
 $g\in L^{2,1}(\R)$,  supp $g\subset B(x_0,r\alpha)$, with  , $0<\alpha<1$, $0<r<\rho/8$.\par
 Then we have
 \begin{eqnarray}\label{eqF2}
 \int_{\R} F(Q,v)\Delta^{-1/4} g dx &\lesssim&\left[\alpha^{1/2}+\left( \frac{r}{\rho}\right)^{1/2}\right]  \|Q\|_{L^2}\|g\|_{L^{2,1}}   \|v\|_{L^{2,\infty}(B(x_0,r))}\\
 &+&\alpha ^{1/2} (\|Q\|_{L^2}+\|Q\|_{L^{\infty}})\|g\|_{L^{2,1}} \sum_{h=1}^{+\infty} 2^{-h/2}\|v\|_{L^{2,\infty}(A_{h,x_0})}\,.\nonumber
 \end{eqnarray}
 \end{Lemma}
 {\bf Proof of Lemma \ref{loc6}.} We just give a sketch of proof.\par
 We write
 \begin{eqnarray*} 
 \int_{\R} F(Q,v)\Delta^{-1/4} g dx&=& \underbrace{\int_{\R} F(Q,\chi_r v)\11_{B(x_0,r/4)}\Delta^{-1/4} g dx}_{(15)}\\
 &+& \underbrace{\int_{\R} F(Q,\chi_r v)\11_{A^\prime_{h,x_0}}\Delta^{-1/4} g dx}_{(16)}\\
 &+& \underbrace{ \int_{\R} F(Q,(1-\chi_r) v)\Delta^{-1/4} g dx}_{(17)}\,.
 \end{eqnarray*}
   To estimate $(15)$ we write $Q=\sum_{h=-2} \tilde\varphi_h Q$ with supp $ \tilde\varphi_h\subseteq B(0,2^{h+1}\rho\setminus B(0,2^{h-1})$ and $ \tilde\varphi _h$ partition of unity.\par
   \begin{eqnarray*}
   (15)&=&\sum_{h=-2}^{\infty}\int_{\R} {\cal{R}}(  \tilde\varphi_h Q) {\cal{R}}(\chi_r v)\11_{B(x_0,r/4)}\Delta^{-1/4} g dx\\[5mm]
   &=&\sum_{h=-2}^{\infty}\int_{\R}{\cal{F}}^{-1}\left[\frac{\cdot}{|\cdot|}\right](\xi)(   \tilde\varphi_h Q)\ast [ {\cal{R}}(\chi_r v)\11_{B(x_0,r/4)}\Delta^{-1/4} g] d\xi\\[5mm]
   &\lesssim&\sum_{h=-2}^{\infty} \|\xi^{-1}\|_{L^\infty(B^c(0,2^h\rho ))}\|   \tilde\varphi_h Q\|_{L^1}\|{\cal{R}}(\chi_r v)\|_{L^1(B(x_0,r/4))}\|\Delta^{-1/4} g\|_{L^\infty(\R)}\\[5mm]
    &\lesssim&\|g\|_{L^{2,1}}\|{\cal{R}}(\chi_r v)\|_{L^{2,\infty}(B(x_0,r/4))} \left( \frac{r}{\rho}\right)^{1/2}\sum_{h=-2}^{\infty}2^{-h/2}\|Q\|_{L^2(A_{h,0})}\\
  &\lesssim& \left(\frac{r}{\rho}\right)^{1/2}\, \|g\|_{L^{2,1}}\|Q\|_{L^2}\| v\|_{L^{2,\infty}(B(x_0,r))}\,.
    \end{eqnarray*}
    Now we write
   \begin{eqnarray*} 
   (16)&=&\sum_{h=-2}^{+\infty}\sum_{k=-2}^{+\infty} \int_{\R}F(\tilde \varphi_k Q,\chi_r v) \11_{A^\prime_{h,x_0}}\Delta^{-1/4} gdx\\[5mm]
   &=&\sum_{h=-2}^{+\infty}\sum_{|k-h|\le 5} \int_{\R}F(\tilde \varphi_ki Q,\chi_r v) \11_{A^\prime_{h,x_0}}\Delta^{-1/4} gdx\\[5mm]
   &+& \sum_{h=-2}^{+\infty}\sum_{|k-h|> 5} \int_{\R}F(\tilde \varphi_k Q,\chi_r v) \11_{A^\prime_{h,x_0}}\Delta^{-1/4} gdx\\[5mm]
   &&\mbox{ by arguing as in $(5)$ and $(6)$}\\[5mm]
    &\lesssim&\|g\|_{L^{2,1}}\|Q\|_{L^2}\left[\left( \frac{r}{\rho}\right)^{1/2}+\alpha^{1/2}\right]\,.
   \end{eqnarray*}
      The estimate of $(16)$ is analogous  to $(2)$ in the proof of Lemma \ref{loc1} and we omit it.~\hfill$\Box$

  \section{ Commutator Estimates}\label{commutators}
 We consider the  Littlewood-Paley   decomposition of unity introduced in Section \ref{defnot}.    
  For every $j\in \Z$ and $f\in{\cal{S}}^{\prime}(\R^n)$  we define the Littlewood-Paley projection operators $P_j$ and $P_{\le j}$ by 
   \begin{eqnarray*}
 \widehat{ P_jf}=\psi_j \hat{f}~~~\widehat{ P_{\le j}f}=\phi_j \hat{f}\,.
 \end{eqnarray*}
 Informally $P_j$ is a frequency projection to the annulus $\{2^{j-1}\le |\xi|\le 2^{j}\}$, while
 $P_{\le j}$ is a frequency projection to the ball $\{|\xi|\le 2^j\}\,.$ We will set
 $f_j=P_j f$ and $f^j=P_{\le j} f$\,.
 
 We observe that  $f^j=\sum_{k=-\infty}^{j} f_k$  and $f=\sum_{k=-\infty}^{+\infty}f_k$ (where the convergence is in ${\cal{S}}^\prime(\R^n)$)\,.\par
 Given $f,g\in {\cal{S}}^\prime(\R)$  we can  split  the product in the following way
 \begin{equation}\label{decompbis}
 fg=\Pi_1(f,g)+\Pi_2(f,g)+\Pi_3 (f,g),\end{equation}
 where
 \begin{eqnarray*}
 \Pi_1(f,g)&=& \sum_{-\infty}^{+\infty} f_j\sum_{k\le j-4} g_k= \sum_{-\infty}^{+\infty} f_j g^{j-4}\,;\\
 \Pi_2(f,g)&=& \sum_{-\infty}^{+\infty} f_j\sum_{k\ge j+4} g_k =\sum_{-\infty}^{+\infty} g_j f^{j-4}\,;\\
 \Pi_3(f,g)&=& \sum_{-\infty}^{+\infty}  f_j\sum_{|k-j|< 4} g_k\,.\
 \end{eqnarray*}
 We observe that   for every $j$ we have 
 $$\mbox{supp${\cal{F}}[f^{j-4}g_j]\subset \{2^{j-2}\le |\xi|\le 2^{j+2}\}$};$$
 $$\mbox{supp${\cal{F}}[\sum_{k=j-3}^{j+3}f_jg_k]\subset \{|\xi|\le 2^{j+5}\}$}\,.$$
 The three pieces of the decomposition \rec{decompbis} are examples of paraproducts. Informally the
 first paraproduct $\Pi_1$ is an operator which allows high frequences of $f$ $(\thicksim 2^j)$ multiplied by  low frequences of $g$ $(\ll 2^j)$ to produce high frequences in the output. The second paraproduct
 $\Pi_2$ multiplies low fequences of $f$ with high frequences of $g$ to produce high fequences in the output. The third paraproduct $\Pi_3$ multiply high frequences of $f$ with high frequences of $g$ to produce comparable or lower frequences in the output. For a presentation of these paraproducts we refer to the reader for instance  to the book \cite{Gra2}\,.
 The following two Lemmae will be often  used   in the sequel. For the proof of the first one   we refer the reader to
 \cite{DLR}\,.
 
 \begin{Lemma}
 For every $f\in {\cal{S}}^{\prime}$ we have
 $$ \sup_{j\in Z}|f^j|\le M(f)\,.$$
 \end{Lemma}
  \begin{Lemma}\label{lemmaprel1}
Let $\psi$ be a Schwartz  radial function such that $supp (\psi)\subset B(0,4)$. Then 
$$
\|\nabla^k {\cal{F}}^{-1}\psi\|_{L^1}\le C_{\psi,n} 4^k\,,
$$
where $C_{\psi,n}$ is a positive constant depending on the $C^2$ norm of $\psi$ and the dimension\,.
\end{Lemma}
{\bf Proof of Lemma \ref{lemmaprel1}.} 
We recall that
$$\nabla^k {\cal{F}}^{-1}\psi (\xi)={\cal{F}}^{-1}[i^kx^k \psi](\xi)\,.$$

We write
$$\int_{\R^n}|\nabla^k {\cal{F}}^{-1}\psi (\xi)|d\xi=\int_{|\xi|\le 1}|\nabla^k {\cal{F}}^{-1}\psi (\xi)|d\xi+\int_{|\xi|\ge 1}|\nabla^k {\cal{F}}^{-1}\psi (\xi)|d\xi\,.$$

 The following estimates hold.\par
 \begin{eqnarray}\label{le1}
\int_{|\xi|\le 1}|\nabla^k {\cal{F}}^{-1}\psi (\xi)|d\xi &\le& \omega_n \|\nabla^k {\cal{F}}^{-1}\psi (\xi)\|_{L^\infty}\\&&\le \omega_n 
\|x^k\psi\|_{L^1}\le\omega_n 4^k\|\psi\|_{L^1}\,,\nonumber
\end{eqnarray}
where $\omega_n=|B_1(0)|\,.$
\begin{eqnarray}\label{ge1}
&&
\int_{|\xi|\ge 1}|\nabla^k {\cal{F}}^{-1}\psi (\xi)|d\xi=\int_{|\xi|\ge 1}(-\frac{1}{|\xi|^2})\left[\int_{\R^n}(\Delta_{x}e^{i\xi x} )\psi(x)(ix)^kdx\right]d\xi\\
&& =\int_{|\xi|\ge 1}(-\frac{1}{|\xi|^2})\left[\int_{\R^n} e^{i\xi x} \Delta_{\xi}(\psi(x)(ix)^k)dx\right]d\xi\nonumber\\
&&
\le \int_{|\xi|\ge 1}\frac{1}{|\xi|^2}d\xi\left(\frac{4^k+2k 4^{k-1}+k(k-1)4^{k-2}}{|x|^2}\|\psi\|_{C^2}\right)\,.\nonumber
\end{eqnarray}
By combining \rec{le1} and \rec{ge1} we obtain
\begin{eqnarray}\label{finalpsi}
&&
\|\nabla^k {\cal{F}}^{-1}\psi (\xi)\|_{L^1} 
\le \omega_n 4^k\|\psi\|_{L^1}\\
&&~~~~+ \frac{4^k+2k 4^{k-1}+k(k-1)4^{k-2}}{|x|^2}\|\psi\|_{C^2}\int_{|\xi|\ge 1}\frac{1}{|\xi|^2}d\xi
\nonumber\\
&&\le C_{\psi,n} 4^k\,. \nonumber
\end{eqnarray}
This concludes the proof of Lemma \ref{lemmaprel1}.~\hfill$\Box$
\par
\medskip
\begin{Lemma}\label{lemmaprel2}
Let $f\in B^0_{\infty,\infty}(\R^n)$. Then for all $k\in\N$  and for all $j\in Z$ we have 
$$
2^{-kj}\|\nabla^k f_j\|_{L^\infty}\le 4^k \|f_j\|_{L^\infty}\,.
$$
\end{Lemma}
{\bf Proof of Lemma \ref{lemmaprel2}.} 
Let $\Psi$ be a Schwartz  radial function such that $\Psi=1$ in $B_2$ and $\Psi=0$ in $B^c(0,4)\,.$\par
Since $supp {\cal{F}}[f_j]\subseteq B_{2^{j+1}}\setminus B_{2^{j-1}}$ we have
\begin{equation}
{\cal{F}}[\nabla^k f_j]\simeq \xi^k {\cal{F}}[f_j]=2^{kj}\psi(2^{-j}\xi)\frac{\xi^k}{2^{k j}}{\cal{F}}[f_j]\,.
\end{equation}
Observe that
\begin{eqnarray*}
&&\|{\cal{F}}^{-1}[\psi(2^{-j}\xi)\frac{\xi^k}{2^{k j}}]\|_{L^1}\\
&&=\|\int_{\R^n}e^{ix\xi}\psi(2^{-j}\xi)\frac{\xi^k}{2^{k j}}d\xi\|_{L^1}\\
&& =2^{nj}\|\int_{\R^n}e^{i2^j x\xi}\psi(\xi){\xi}^kd\xi\|_{L^1}\\
&&=2^{nj}\|\nabla^k{\cal{F}}^{-1}[\psi](2^j \cdot)\|_{L^1}\,.
\end{eqnarray*}
Thus
\begin{eqnarray*}
&&
 2^{-kj}\|\nabla^k f_j\|_{L^\infty}\le \|{\cal{F}}^{-1}[\psi(2^{-j}\xi)\frac{\xi^k}{2^{k j}}]\ast f_j\|_{L^\infty}\\
&&\le
 \|{\cal{F}}^{-1}[\psi(2^{-j}\xi)\frac{\xi^k}{2^{k j}}]\|_{L^1}\|f_j\|_{L^\infty}\\
&&
= \|{\cal{F}}^{-1}[\nabla ^k \psi]\|_{L^1}\|f_j\|_{L^\infty}
\le C_{\psi} 4^k\|f_j\|_{L^\infty}\,.~~\hfill\Box
\end{eqnarray*}

 %%%%%%%%%%%%%%%%%%%%%%%%%%
Now we start with a series of preliminary Lemmae which will be crucial in the construction of the gauge
$P$ in Section \ref{constrP} 
\begin{Lemma}\label{lemmaprel1bis}
Let $a\in \dot W^{1/2,r}(\R)$, $r<2$ and $b\in{ \dot {H}}^{1/2}(\R)$. Then
\begin{eqnarray}\label{comm1}
&&\|\Delta^{1/4} (ab)-a\Delta^{1/4}b-(\Delta^{1/4} a) b\|_{L^r(\R)}\le C \|a\|_{\dot W^{1/2,r}(\R)}\|b\|_{\dot H^{1/2}(\R)}\,.\nonumber\end{eqnarray}
\end{Lemma}
{\bf\noindent Proof of Lemma \ref{lemmaprel1bis}.}\par
{\bf $\bullet$ Estimate of $\|\Pi_2(\Delta^{1/4} (ab))\|_{L^r}$\,.} 
\begin{eqnarray*}
&&
\|\sum_j\Delta^{1/4}(a^{j-4}b_j)\|^r_{L^r} \lesssim \int_{\R}\left(\sum_j 2^j |a^{j-4}|^2|b_j|^2\right)^{r/2} dx\\
&& \lesssim \int_{\R} \sup_j|a^{j-4}|^r\left(\sum_j2^j|b_j|^2\right)^{r/2}\\[5mm]
&&\mbox{by  H\"older Inequality}\\[5mm]
&& \lesssim\left( \int_{\R} \sup_j|a^{j-4}|^{\frac{2r}{2-r}}\right)^{\frac{2-r}{2}}
\left( \int_{\R} \sum_j|b_j|^2 dx\right)^{r/2}\\[5mm]
&& \lesssim \|a\|^r_{L^{{\frac{2r}{2-r}}}}\|b\|^r_{\dot H^{1/2}} \lesssim \|a\|^r_{\dot W^{1/2,r}}\|b\|^r_{\dot H^{1/2}}\,.
\end{eqnarray*}
In the last inequality we use the embedding $\dot W^{1/2,r}(\R)\hookrightarrow L^{{\frac{2r}{2-r}}}(\R)\,,$(see for instance \cite{AD2}).\par
{\bf $\bullet$ Estimate of $\|\Pi_2(a\Delta^{1/4}b)\|_{L^r}$} \,.
\begin{eqnarray*}
&&
\|\Pi_2(a\Delta^{1/4}b)\|_{L^r}^r \lesssim \int_{\R}\left(\sum_j |a^{j-4}|^2|\Delta^{1/4} b_j|^2 dx\right)^{r/2}dx \\
&& \lesssim  \int_{\R}\sup_j|a^{j-4}|^2\left(\sum_j |\Delta^{1/4} b_j|^2\right)^{r/2}dx \\
&& \lesssim \|a\|^{r}_{\dot W^{1/2,r}} \|b\|^r_{\dot H^{1/2}}\,.
\end{eqnarray*}
{\bf $\bullet$ Estimate of $\|\Pi_2((\Delta^{1/4} a) b)\|_{L^r}$} \,.
\begin{eqnarray*}
&& \|\Pi_2((\Delta^{1/4} a) b)\|^r_{L^r}\lesssim 
 \int_{\R}\left(\sum_j |\Delta^{1/4}a^{j-4}|^2|b_j|^2\right)^{r/2} dx\\
 &&\lesssim 
 \int_{\R} \sup_j \left(2^{-j/2}|\Delta^{1/4} a^{j-4}|\right)^r\left(\sum_j 2^j|b_j|^2\right)^{r/2} dx \\
 &&\lesssim 
 \left( \int_{\R} \sup_j \left(2^{-j/2}|\Delta^{1/4} a^{j-4}|\right)^{\frac{2r}{2-r}}\right)^{\frac{2-r}{2}}
 \left(\int_{\R} \sum_j 2^j |b_j|^2 dx\right)^{r/2}\\[5mm]
 &&\lesssim 
 \left [\int_{\R}\left(\sum_{j} 2^{-j}|\Delta^{1/4} a^{j-4}|^2\right)^{\frac{r}{2-r}}dx\right]^{\frac{2-r}{r}
 }\|b\|^r_{\dot H^{1/2}}\\[5mm]
  &&\lesssim \|a\|_{\dot W^{1/2,r}}\|b\|^r_{\dot H^{1/2}}.
  \end{eqnarray*}
{\bf $\bullet$ Estimate of $\|\Pi_3(\Delta^{1/4} (ab))\|_{L^r}$}\,. 
\begin{eqnarray}\label{p3ab1}
&&
\|\Pi_3(\Delta^{1/4} (ab))\|_{L^r}^r\simeq\sup_{\|h\|_{L^{r^\prime}}\le 1}\int_{\R}\sum_j( \Delta^{1/4} h) a_jb_j dx\nonumber\\
&&
=\sup_{\|h\|_{L^{r^\prime}}\le 1}\left[\int_{\R}\sum_{j}\sum_{|k-j|\le 3}( \Delta^{1/4} h_k) a_jb_j dx+
\int_{\R}\sum_{j}( \Delta^{1/4} h^{j-4}) a_jb_j dx\right]\,.
\end{eqnarray}
Now we estimate the last two terms in \rec{p3ab1}.
\begin{eqnarray*}
&&
\int_{\R}\sum_{j} \Delta^{1/4} h^{j-4} a_jb_j dx\lesssim 
 \int_{\R}\left(\sum_j 2^{-j} |\Delta^{1/4} h^{j-4} |^2\right)^{1/2}\left(\sum_j 2^j|a_j|^2|b_j|^2\right)^{1/2} dx\\[5mm]
 &&\lesssim 
 \left[ \int_{\R}\left(\sum_j 2^{-j} |\Delta^{1/4} h^{j-4} |^2\right)^{r^{\prime}/2}dx\right]^{1/r^{\prime}}\left[\int_{\R}  \left(\sum_j 2^j|a_j|^2|b_j|^2\right)^{r/2}dx\right]^{1/r}\\[5mm]
 &&\lesssim 
 \|h\|_{L^{r^\prime}}\|b\|_{B^0_{\infty,\infty}}\|a\|_{\dot W^{1/2,r}}\,.
 \end{eqnarray*}
 The estimate of $\int_{\R}\sum_{j}\sum_{|k-j|\le 3}( \Delta^{1/4} h_k) a_jb_j dx$ is similar and we omit it.\par
 \medskip
 {\bf $\bullet$ Estimate of $\|\Pi_3(a\Delta^{1/4} b)\|_{L^r}$}\,.
\begin{eqnarray*}
&&
\|\Pi_3(a\Delta^{1/4} b)\|_{L^r}^r 
\lesssim
\int_{\R}\left|\sum_j a_j\Delta^{1/4} b_j\right|^r dx\\[5mm]
&&\lesssim
\int_{\R}\left(\sum_j a_j^2\right)^{r/2}\left(\sum_j| \Delta^{1/4}b_j|^2\right)^{r/2}dx\\[5mm]
&&\mbox{by Cauchy-Schwartz Inequality}\\[5mm]
&&\lesssim
\left(\int_{\R}\left(\sum_j a_j^2\right)^{r/(2-r)}dx\right)^{(2-r)/2}
\left(\int_{\R}\left(\sum_j |\Delta^{1/4}b_j|^2\right)^{r/(2-r)}dx \right)^{(2-r)/2}\\[5mm]
&&\lesssim
\|a\|^r_{L^{\frac{2r}{2-r}}}\|b\|^r_{\dot H^{1/2}}\lesssim \|a\|^r_{\dot W^{1/2,r}}\|b\|^r_{\dot H^{1/2}}.
  \end{eqnarray*}
 {\bf $\bullet$ Estimate of $\|\Pi_3((\Delta^{1/4} a) b)\|_{L^r}$}\,. 
 \begin{eqnarray*}
&&
\|\Pi_3((\Delta^{1/4} a) b)\|_{L^r}^r \lesssim\int_{\R} \left(\sum_j 2^{-j}|\Delta^{1/4} a_i|^2\right)^{r/2}\left(\sum_{j} 2^j|b_j|^2\right)^{r/2}dx\\
&&\lesssim
\|a\|_{\dot W^{1/2,r}}\|b\|^r_{\dot H^{1/2}}.
  \end{eqnarray*} 
   {\bf $\bullet$ Estimate of $\|\Pi_1(a\Delta^{1/4} b)\|_{L^r}$}\,. 
 \begin{eqnarray*}
&&
\|\Pi_1(a\Delta^{1/4} b)\|_{L^r}\simeq \sup_{\|h\|_{L^{r^\prime}}\le 1}\int_{\R}\sum_j \Delta^{1/4} b^{j-4} a_jh_j dx\\[5mm]
&&\lesssim
\int_{\R}\sup_j|\Delta^{1/4} b^{j-4}|\left(\sum_{j} a_j^2\right)^{1/2}\left(\sum_{j} h_j^2\right)^{1/2}\\[5mm]
&&\mbox{by generalized H\"older Inequality: $\frac{1}{2}+\frac{1}{r^{\prime}}+\frac{2-r}{2r}=1$}\\[5mm]
&&
\lesssim
\|b\|_{\dot H^{1/2}} \|a\|_{L^{{\frac{2r}{2-r}}}}\|h\|_{L^{r^\prime}}
\lesssim
\|b\|_{\dot H^{1/2}}  \|a\|_{\dot W^{1/2,r}}\,.
    \end{eqnarray*} 
      {\bf $\bullet$ Estimate of $\|\Pi_1(\Delta^{1/4} (ab)-(\Delta^{1/4} a)b)\|_{L^r}$\,.}
       \begin{eqnarray}\label{piab}
       &&
\|\Pi_1(\Delta^{1/4} (ab)-(\Delta^{1/4} a)b)\|_{L^r} \\
&&=\sup_{\|h\|_{L^{r^\prime}}\le 1}\int_{\R}\sum_j \sum_{|k-j\le 3}|h_k(\Delta^{1/4}(a_j b^{j-4})-\Delta^{1/4} a_j b^{j-4}) dx\nonumber\\
&&=\sup_{\|h\|_{L^{r^\prime}}\le 1}\int_{\R}\sum_j  \sum_{|k-j\le 3}b^{j-4}(\Delta^{1/4}(h_k) a_{j})  -\Delta^{1/4} (a_j) h_{k}) d\xi \nonumber 
 \end{eqnarray}
  \begin{eqnarray*}
&&\lesssim
\sup_{\|h\|_{L^{r^\prime}}\le 1} \int_{\R}\sum_j\sum _{|k-j|\le 3} {\cal{F}}[b^{j-4}]
{\cal{F}}[\Delta ^{1/4}h_k  a_j- \Delta^{1/4} a_j h_{k}) d\xi \nonumber\\ && 
 =
\sup_{\|h\|_{L^{r^\prime}}\le 1} \int_{\R^n}\sum_j\sum _{|k-j|\le 3} {\cal{F}}[b^{j-4}](\xi)\nonumber\\
 &&\left[\int_{\R^n}({\cal{F}}[h_j](\eta){\cal{F}}[a_j](\xi-\eta)(|\eta|^{1/2}-|\xi-\eta|^{1/2})d\eta\right]d\xi\,.\nonumber
 \end{eqnarray*}
 Now we observe that in \rec{piab} we have $|\xi|\le 2^{j-3}$ and $2^{j-2}\le |\eta|\le 2^{j+2}$.
Thus $|\displaystyle\frac{\xi}{\eta}|\le \frac{1}{2}\,.$ 
Hence
\begin{eqnarray}\label{estkern}
|\eta|^{1/2}-|\xi-\eta|^{1/2}&=&|\eta|^{1/2}[1-|1-\frac{\xi}{\eta}|^{1/2}]\\
&=&|\eta|^{1/2}\frac{\xi}{\eta}[1+|1-\frac{\xi}{\eta}|^{1/2}]^{-1}\nonumber\\
&=&|\eta|^{1/2}\sum_{k=0}^\infty\frac{c_k}{k!}(\frac{\xi}{\eta})^{k+1}\,.\nonumber
\end{eqnarray}
 We may suppose that $\sum_{k=0}^\infty\frac{c_k}{k!}(\frac{\xi}{\eta})^{k+1}$ is convergent if $|\displaystyle\frac{\xi}{\eta}|\le \frac{1}{2}\,,$ otherwise one may consider a
 different Littlewood-Paley decomposition by replacing the exponent $j-4$ with $j-s$, $s>0$ large enough. 
 We    introduce the following notation: for every $k\ge 0$ we set
$$S_k g={\cal{F}}^{-1}[\xi^{-(k+1)}|\xi|^{1/2} {\cal{F}} g].$$
We note that  if $h\in B^s_{\infty,\infty}$ then
$S_k h\in B^{s+1/2+k}_{\infty,\infty}$ and if $h\in H^{s}$ then $S_k h\in H^{s+1/2+k}\,.$
Moreover if $Q\in H^{1/2}$ then $\nabla^{k+1}(Q)\in H^{-k-1/2}\,.$\par
Next we continue with the proof of \rec{piab}\,.
 \begin{eqnarray}\label{pi1crochet1}
&&
 \sup_{\|h\|_{L^{r^\prime}}\le 1}\int_{\R}\sum_j \sum_{|k-j\le 3}h_k(\Delta^{1/4}(a_j b^{j-4})-(\Delta^{1/4} a_j )b^{j-4}) dx\nonumber\\[5mm]
 &&\lesssim
   \sup_{\|h\|_{L^{r^\prime}}\le 1} 
\sum_{\ell=0}^\infty\frac{c_{\ell}}{\ell!}\int_{\R}
\sum_j\sum _{|k-j|\le 3}[\nabla^{\ell+1}b^{j-4} [(S_{\ell}h_k) a_j)]( x) d x\nonumber\\[5mm]
&&\mbox{by  Lemma \ref{lemmaprel2} }\nonumber\\[5mm]
 &&\lesssim
 \sup_{\|h\|_{L^{r^\prime}}\le 1} 
\sum_{\ell=0}^\infty\frac{c_{\ell}}{\ell!}2^{-4\ell}4^{\ell+1}\|b\|_{B_{\infty,\infty}^0}
\int_{\R}
\sum_j\sum _{|k-j|\le 3}2^{(\ell+1)j} |[(S_{\ell}h_k)a_j)]( x) |d x\nonumber\\[5mm]
&&\lesssim
 \sup_{\|h\|_{L^{r^\prime}}\le 1} 
\sum_{\ell=0}^\infty\frac{c_{\ell}}{\ell!}2^{-4\ell}4^{\ell+1}\|b\|_{B_{\infty,\infty}^0}
\int_{\R}\sum_j |(2^{(1/2+\ell)j}S_{\ell}h_j)|\,|( 2^{j/2} a_j) |dx\\[5mm]
&&\mbox{by Schwartz Inequality}\nonumber\\[5mm]
&&\lesssim
 \sup_{\|h\|_{L^{r^\prime}}\le 1} 
\sum_{\ell=0}^\infty\frac{c_{\ell}}{\ell!}2^{-4\ell}4^{\ell+1}\|b\|_{B_{\infty,\infty}^0}
\int_{\R}(\sum_j 2^{2j(1/2+\ell)}|S_{\ell}h_j|^2)^{1/2}  (\sum_j 2^ja^2_j )^{1/2}dx\nonumber\\[5mm]
&&\mbox{by H\"older  Inequality\nonumber}\\[5mm]
&&\lesssim
 \sup_{\|h\|_{L^{r^\prime}}\le 1} 
\sum_{\ell=0}^\infty\frac{c_{\ell}}{\ell!}2^{-4\ell}4^{\ell+1}\|b\|_{B_{\infty,\infty}^0}\nonumber\\[5mm]&&~~~
\left(\int_{\R}(\sum_j 2^{2j(1/2+\ell)}|S_{\ell}h_j|^2)^{r^{\prime}/2}\right)^{1/r^{\prime}}
\left(\int_{\R}(\sum_j  2^{j} a_j^2)^{r/2}\right)^{1/r}\nonumber\\[5mm]
&&\lesssim
 \sup_{\|h\|_{L^{r^\prime}}\le 1} \sum_{\ell=0}^\infty\frac{c_{\ell}}{\ell!}2^{-2\ell} \|b\|_{B_{\infty,\infty}^0}\|h\|_{L^{r^\prime}}\|\Delta^{1/4}a\|_{L^r}\nonumber\\[5mm]
 &&\lesssim
 \|b\|_{\dot H^{1/2}}  \|a\|_{\dot W^{1/2,r}}\, \nonumber
    \end{eqnarray}\par
    This concludes the proof of Lemma \ref{lemmaprel1bis}\,.~~~~~~\hfill$\Box$
   \begin{Lemma}\label{lemmaprel2bis}
Let $1<r<2$, $a\in \dot W^{1/2,r}(\R)$ and $b\in{ \dot {H}}^{1/2}(\R)\cap L^{\infty}(\R)$\,.
Then
$$
\|ab\|_{\dot W^{1/2,r}}\le C\|a\|_{\dot W^{1/2,r}}(\|b\|_{\dot H^{1/2}}+\|b\|_{L^{\infty}})\,.$$
\end{Lemma}
{\bf Proof of Lemma \ref{lemmaprel2bis}\,.}
 {\bf $\bullet$ Estimate of $\|\Pi_1(\Delta^{1/4}(ab))\|_{L^r}\,.$}\par
 \begin{eqnarray*}
 &&
 \|\sum_j\Delta^{1/4} (a_j b^{j-4})\|^r_{L^r}
 \lesssim
 \int_{\R}(\sum_j|a_j|^2|b^{j-4}|^2)^{r/2}dx \\
 &&
~~~~\int_{\R} \sup_j |b^{j-4}|^r \left(\sum_j 2^j |a_j|^2\right)^{r/2}dx\\
  &&
 \lesssim
 \int_{\R} |M(b)|^r \left(\sum_j 2^j |a_j|^2\right)^{r/2} dx
 \le \|b\|^r_{L^{\infty}}
 \int_{\R} \left(\sum_j 2^j |a_j|^2\right )^{r/2} dx\\
   &&
 \lesssim
\|b\|^r_{L^{\infty}} \|a\|^r_{\dot W^{1/2,r}}\,.
\end{eqnarray*}
\par
 {\bf $\bullet$ Estimate of $\|\Pi_2\Delta^{1/4}(ab)\|_{L^r}\,.$}\par
 \begin{eqnarray*}
 &&
 \|\sum_j\Delta^{1/4}( a^{j-4} b_j)\|_{L^r}
 \simeq \sup_{\|h\|_{L^{r^\prime}}\le 1} 
 \int_{\R}\sum_j  a^{j-1}b_j \Delta^{1/4} h_j\\
   &&
 \lesssim
\sup_{\|h\|_{L^{r^\prime}}\le 1} 
\int_{\R}\sup_j |a^{j-4}|\left(\sum_j 2^j|b_j|^2\right)^{1/2}
\left(\sum_j|h_j|\right)^{1/2} dx\\[5mm]
&&
 \lesssim
\sup_{\|h\|_{L^{r^\prime}}\le 1} 
\int_{\R} |M(a)|\left(\sum_j 2^j|b_j|^2\right)^{1/2}
\left(\sum_j|h_j|\right)^{1/2} dx\\[5mm]
&&\mbox{by generalized H\"older Inequality: $\frac{1}{r^{\prime}}+\frac{1}{2}+\frac{2-r}{2r}=1$}\\[5mm]
&&
 \lesssim
 \|b\|_{\dot H^{1/2}} \|a\|_{\dot W^{1/2,r}}\,.
 \end{eqnarray*}
\par
 {\bf $\bullet$ Estimate of $\|\Pi_3(\Delta^{1/4}(ab)\|_{L^r}\,.$}\par
 \begin{eqnarray*}
 &&
 \|\sum_j\Delta^{1/4}(a_jb_j)\|_{L^r}\\
 &&\simeq \sup_{\|h\|_{L^{r^\prime}}\le 1}\left[\int_{\R}\sum_j\sum_{|k-j|\le 3} \Delta^{1/4}(a_jb_j) h_k dx+
 \int_{\R}\sum_j \Delta^{1/4}(a_jb_j) h^{j-4} dx\right]\\
 &&
 =
 \sup_{\|h\|_{L^{r^\prime}}\le 1}\left[\int_{\R}\sum_j\sum_{|k-j|\le 3} (a_jb_j) \Delta^{1/4}h_k dx+
 \int_{\R}\sum_j (a_jb_j)\Delta^{1/4} h^{j-4} dx\right]
   \end{eqnarray*}
   We estimate the term $\int_{\R}\sum_j (a_jb_j)\Delta^{1/4} h^{j-4} dx\,.
$
\begin{eqnarray*}
&&
|\int_{\R}\sum_j (a_jb_j)\Delta^{1/4} h^{j-4} dx|\\[5mm]
&&
 \lesssim \|b\|_{B^0_{\infty,\infty}}\int_{\R}\left(\sum_{j} 2^{-j}|\Delta^{1/4} h^{j-4}|^2\right)^{1/2}
 \left(\sum_{j} 2^{j}a_j^2\right)^{1/2}dx\\[5mm]
  && \lesssim \|b\|_{B^0_{\infty,\infty}}\left(\int_{\R}\left(\sum_{j} 2^{-j}|\Delta^{1/4} h^{j-4}|^2\right)^{r^{\prime}/2}dx \right)^{1/r^{\prime}}\left(\int_{\R} \left(\sum_{j} 2^{j}a_j^2\right)^{r/2}dx\right)^{1/r}\\ [5mm]
&&
 \lesssim \|b\|_{B^0_{\infty,\infty}}\|h\|_{L^{r^{\prime}}} \|a\|_{\dot W^{1/2,r}}\,.
 \end{eqnarray*}
 The term $\int_{\R}\sum_j\sum_{|k-j|\le 3} (a_jb_j) \Delta^{1/4}h_k dx$ is estimated in a similar way.
 Thus
  we get
   \begin{eqnarray*}
 &&
 \|\sum_j\Delta^{1/4}(a_jb_j)\|_{L^r}\lesssim \|b\|_{\dot H^{1/2}} \|a\|_{\dot W^{1/2,r}}\,. 
 \end{eqnarray*} 
  This concludes the proof of Lemma \ref{lemmaprel2bis}\,.~\hfill $\Box$
  \medskip
   \begin{Lemma}\label{lemmaprel3}
Let $1<r<2<q$, $a\in \dot W^{1/2,r}(\R)$ and $b\in \dot W^{1/2,q}(\R)$ and $t=\frac{2rq}{2r+q(2-r)}\,.$ Then
\begin{eqnarray}\label{comm3}
&&\|\Delta^{1/4} (ab) -(\Delta^{1/4} a) b\|_{L^t(\R)}\le C \|a\|_{\dot W^{1/2,r}(\R^n)}\|b\|_{\dot W^{1/2,q}(\R)}\,.\nonumber\end{eqnarray}
\end{Lemma}
\par
{\bf Proof of Lemma \ref{lemmaprel3}\,.}\par
{\bf $\bullet$ Estimate of $\|\Pi_2(\Delta^{1/4}(ab))\|_{L^t}\,.$}\par
\begin{eqnarray*}
&&\|\sum_j \Delta^{1/4}(a^{j-4}b_j)\|^t_{L^t}
\lesssim \int_{\R} \left(\sum_j 2^j |a^{j-4}|^2|b_j|^2\right)^{t/2} dx\\[5mm]
&&
\lesssim
\int_{\R}\sup_j|a^{j-4}|^t\left(\sum_{j} 2^j|b_j|^2\right)^{t/2}\\[5mm]
&&
\lesssim
\left(\int_{\R}M(a)^{\frac{tq}{q-t}}dx\right)^{1-\frac{t}{q}}\left(\int_{\R}(\sum_j 2^j|b_j|^2)^{q/2} dx\right)^{t/q}\\[5mm]
&&
\lesssim
\|a\|^t_{L^{\frac{2r}{2-r}} }\|b\|^t_{\dot W^{1/2,q}}\,.
\end{eqnarray*}
In the above expression we use the fact that $\frac{tq}{q-t}=\frac{2r}{2-r}\,.$\par
{\bf $\bullet$ Estimate of $\|\Pi_2((\Delta^{1/4}a) b)\|_{L^t}\,.$}\par
\begin{eqnarray*}
&&\|\sum_j( \Delta^{1/4}a^{j-4})b_j\|^t_{L^t}\\
 &&
\lesssim
 \int_{\R} \left(\sup_j 2^{-j/2}|\Delta^{1/4} a^{j-4}|\right)^t \left(\sum_j 2^j |b_j|^2\right)^{t/2} dx\\[5mm]
 && \lesssim
  \int_{\R} \left(\sum_j 2^{-j}|\Delta^{1/4} a^{j-4}|^2\right)^{t/2} \left(\sum_j 2^j |b_j|^2\right)^{t/2} dx\\[5mm]
 && \lesssim
\left( \int_{\R} \left(\sum_j 2^{-j}|\Delta^{1/4} a^{j-4}|^2\right)^{tq/2(q-t)}\right)^{1-t/q}
\left(\left(\sum_j 2^j |b_j|^2\right)^{q/2} dx\right)^{t/q}\\[5mm]
&&
\lesssim \|a\|^t_{L^{tq/q-t} }\|b\|^t_{\dot W^{1/2,q}}
\lesssim\|a\|^t_{\dot W^{1/2,r}} \|b\|^t_{\dot W^{1/2,q}}\,.
\end{eqnarray*}
\par
{\bf $\bullet$ Estimate of $\|\Pi_3(\Delta^{1/4}(ab))\|_{L^t}\,.$}\par
\begin{eqnarray*}
&&\|\sum_j \Delta^{1/4}(a_jb_j)\|_{L^t}
\simeq\sup_{\|h\|_{L^{t^\prime}}\le 1}\int_{\R} \Delta^{1/4} h\sum_{j} a_j b_j dx\\
&& 
\lesssim \sup_{\|h\|_{L^{t^\prime}}\le 1}\left[\int_{\R} \sum_j\sum_{|j-k|\le 4 }\Delta^{1/4} h_k a_j b_j dx+
\int_{\R} \sum_j \Delta^{1/4} h^{j-4}a_j b_j dx\right]\,.
\end{eqnarray*}
We estimate the term  $\int_{\R} \sum_j \Delta^{1/4} h^{j-4}a_j b_j dx$.
\begin{eqnarray*}
&&
\int_{\R} \sum_j \Delta^{1/4} h^{j-4}a_j b_j dx\\[5mm]
&& \lesssim
\int_{\R} \sup_{j} \left(2^{-j/2} |\Delta^{1/4} h^{j-4}\right)|\sum_{j} 2^{j/2} |a_j\|b_j| dx\\[5mm]
&&\lesssim
\int_{\R}\left(\sum_{j} |\Delta^{1/4} h^{j-4}|^2\right)^{1/2}\left(\sum_j |a_j|^2\right)^{1/2}\left(\sum_j 2^{j}|b_j|^2\right)^{1/2} dx\\[5mm]
&&
\lesssim
\left[\int_{\R}\left(\sum_{j} 2^{-j}|\Delta^{1/4} h^{j-4}|^2\right)^{t^{\prime}/2}\right]^{1/t^{\prime}}
\left[\int_{\R} \left(\sum_j |a_j|^2\right)^{t/2}\left(\sum_j 2^{j}|b_j|^2\right)^{t/2} \right]^{1/t}\\[5mm]
&&
\lesssim
\|h\|_{L^{t^{\prime}}\e 1}
\left[\int_{\R} \left(\sum_j |a_j|^2\right)^{tq/2(q-t)}dx\right]^{\frac{q-t}{qt}}
\left[\int_{\R} \left(\sum_j 2^j |b_j|^2\right)^{q/2}dx\right]^{1/q}\\
&&
\lesssim \|a\|_{tq/q-t} \|b\|_{W^{1/2,q}}\\[5mm]
&&
\lesssim \|a\|_{W^{1/2,r}} \|b\|_{W^{1/2,q}}\,.
\end{eqnarray*}
The estimate of $\int_{\R} \sum_j\sum_{|j-k|\le 4 }\Delta^{1/4} h_k a_j b_j dx$ is similar.\par
{\bf $\bullet$ Estimate of $\|\Pi_3((\Delta^{1/4} a)b)\|_{L^t}\,.$}\par
\begin{eqnarray*}
&&\|\sum_j (\Delta^{1/4}a_j)b_j)\|^t_{L^t}
\lesssim\int_{\R}|\sum_j\Delta^{1/4} a_j b_j|^t
\\[5mm]
&&
\lesssim
\int_{\R}\left(\sum_j 2^{-j}|\Delta^{1/4} a_j|^2\right)^{t/2}\left(\sum_j 2^{j}b_j|^2\right)^{t/2} dx\\[5mm]
&&
\lesssim \|a\|^t_{tq/q-t} \|b\|^t_{W^{1/2,q}}
\lesssim \|a\|^t_{W^{1/2,r}} \|b\|^t_{W^{1/2,q}}\,.
\end{eqnarray*}

{\bf $\bullet$ Estimate of $\|\Pi_2(\Delta^{1/4}( ab)-(\Delta^{1/4} a) b)\|_{L^t}\,.$}\par
\begin{eqnarray}\label{provvbis}
&&
\|\sum_j(\Delta^{1/4}( ab)-(\Delta^{1/4} a) b)\|\|_{L^t}\\[5mm]
&&
=\sup_{\|h\|_{L^{t^{\prime}}}\le 1}
\int_{\R}\sum_j h_j [\Delta^{1/4} (a_jb^{j-4})-(\Delta^{1/4} a_j) b^{j_4}] dx\nonumber\\[5mm]
&&
=
\sup_{\|h\|_{L^{t^{\prime}}}\le 1}
\int_{\R}\sum_j b^{j-4} [(\Delta^{1/4} h_j)a_j -h_j(\Delta^{1/4} a_j) ] dx\nonumber\\[5mm]
&&
=
\sup_{\|h\|_{L^{t^{\prime}}}\le 1}\int_{\R}\sum_j {\cal{F}}[b]^{j-4}(\eta)
(\int_{\R} {\cal{F}}[h]_j(\xi){\cal{F}}[a]_j(\eta-\xi)[|\xi|^{1/2}-|\eta-\xi|^{1/2}] d\xi)d\eta\,.\nonumber
\end{eqnarray}
Now we argue as  in \rec{pi1crochet1}
\begin{eqnarray*}
&&\rec{provvbis}
\lesssim
   \sup_{\|h\|_{L^{t^\prime}}\le 1} 
\sum_{\ell=0}^\infty\frac{c_{\ell}}{\ell!}\int_{\R}
\sum_j\sum _{|k-j|\le 3}[\nabla^{\ell+1}b^{j-4} [(S_{\ell}h_k) a_j)]( x) d x\nonumber\\[5mm]
&&
 \lesssim
   \sup_{\|h\|_{L^{t^\prime}}\le 1} 
\sum_{\ell=0}^\infty\frac{c_{\ell}}{\ell!}\int_{\R}
\sum_j[2^{-(\ell+1/2)j}\nabla^{\ell+1}b^{j-4} [2^{j(\ell+1/2)}(S_{\ell}h_j) a_j)]( x) d x\nonumber\\[5mm]
  &&\lesssim
 \sup_{\|h\|_{L^{t^\prime}}\le 1} 
\sum_{\ell=0}^\infty\frac{c_{\ell}}{\ell!}\int_{\R}\sup_j [2^{j(\ell+1/2)}(S_{\ell}h_j) ]\\[5mm]
&&
~~~\left(\sum_j |a_j|^2\right)^{1/2}\left(\sum_j 2^{-2(\ell+1/2)j}|\nabla^{\ell+1}b^{j-4}|^2\right)^{1/2} dx
\end{eqnarray*}
\begin{eqnarray*}
&&\lesssim
\sup_{\|h\|_{L^{t^{\prime}}}\le 1}\sum_{\ell=0}^\infty\frac{c_{\ell}}{\ell!}2^{-2\ell}
\int_{\R} \left(\sum_j 2^{-2(\ell+1/2)j} |S_{\ell}h_j|^2\right)^{1/2}\\[5mm]
&&~~~\left(\sum_j |a_j|^2\right)^{1/2}\left(\sum_j2^{-2(\ell+1/2)j}|\nabla^{\ell+1} b^{j-4}|^2\right)^{1/2} dx\\[5mm]
&& \lesssim
\sup_{\|h\|_{L^{t^{\prime}}}\le 1}
(\int_{\R} (\sum_j 2^{j} |\Delta^{-1/4} h_j|^2)^{t^{\prime}/2})^{t^{\prime}}\\[5mm]
&&~~~
\left[\int_{\R} \left(\sum_j |a_j|^2\right)^{qt/2(q-t)}\right]^{q-t/qt}
\left[\int_{\R} \left(\sum_j |b_j|^2\right)^{q/2}\right]^{1/q}\\[5mm]
&& \lesssim
\|a\|_{L^{qt/q-t}}\|b\|_{W^{1/2,q}}
\lesssim \|a\|_{W^{1/2,r}} \|b\|_{W^{1/2,q}}\,. 
\end{eqnarray*}
This concludes the proof of Lemma \ref{lemmaprel3}.~~\hfill$\Box$
%%%%%%%%%%%% QUI
\begin{Lemma}\label{lemmaprel4}
Let $a\in \dot H^{1/2}(\R^n)\cap L^\infty(\R^n)$,   $b\in W^{1/2,q}(\R^n)$, $2<q<+\infty\,.$ Then
\begin{eqnarray}\label{comm4}
&&\|\Delta^{1/4} (ab)-(\Delta^{1/4} a) b\|_{L^q(\R^n)} \le \|b\|_{\dot W^{1/2,q}(\R^n)}\left[\|a\|_{\dot H^{1/2}(\R^n)}+\|a\|_{L^{\infty}(\R^n)}\right]\,.\nonumber\end{eqnarray}
\end{Lemma}
{\bf Proof of Lemma \ref{lemmaprel4}.}\par
\noindent{\bf $\bullet$ Estimate of $\|\Pi_1(\Delta^{1/4}(ab))\|^q_{L^q}\,.$} 
\begin{eqnarray*}
\|\sum_j 
\Delta^{1/4}(a^{j-4}b_j)\|^q_{L^q}& \lesssim&
\int_{\R}\left(\sum_j 2^j |a^{j-4}|^2|b_j|^2\right)^{q/2}\\[5mm]
& \lesssim&
\|a\|_{L^{\infty}}^q\|b\|^q_{W^{1/2,q}}\,.
\end{eqnarray*}
{\bf $\bullet$ Estimate of $\|\Pi_1((\Delta^{1/4}a )b)\|^q_{L^q}\,.$} 
\begin{eqnarray*}
&&
\|\sum_j \Delta^{1/4} a^{j-4} b_j\|_{L^q}^q
\lesssim 
\int_{\R} \left(\sum_j |\Delta^{1/4} a^{j-4}|^2|b_j|^2\right)^{q/2}\\[5mm]
&&
\lesssim
\|\sup_j 2^{-j/2}|\Delta^{1/4} a^{j-4}|^q_{L^\infty}\int_{\R}\left(\sum_j 2^j |b_j|^2\right)^{q/2} dx\\[5mm]
&&
\lesssim 
\|b\|_{W^{1/2,q}}^q\|a\|_{B^0_{\infty,\infty}}\,.
\end{eqnarray*}

{\bf $\bullet$ Estimate of $\|\Pi_3(\Delta^{1/4}(a b))\|_{L^q}\,.$}\par

\begin{eqnarray*}
&&
\|\sum_j 
\Delta^{1/4}(a_jb_j)\|_{L^q} =\sup_{\|h\|_{L^{q^{\prime}}}\le 1}\int_{\R}(\Delta^{1/4} h) \sum_j a_j b_j dx\\[5mm]
&&
=\sup_{\|h\|_{L^{q^{\prime}}}\le 1}\left[\int_{\R}  \sum_j \sum_{|k-j|\le 4 } (\Delta^{1/4} h_k) a_j b_j dx+\int_{\R}  \sum_j  (\Delta^{1/4} h^{j-4}) a_j b_j dx\right]\,.
\end{eqnarray*}
We estimate the last term:
\begin{eqnarray*}
&&
\int_{\R}  \sum_j \sum_{|k-j|\le 4 }\Delta^{1/4} h^{j-4} a_j b_j dx\\[5mm]
&&
\lesssim \|a\|_{B^0_{\infty,\infty}}\int_{\R}\left|\sum_j 2^{-j}|\Delta^{1/4} h^{j-4}|^2\right|^{1/2}\int_{\R}\left|\sum_j 2^{j}|b_j|^2\right|^{1/2}\\[5mm]
&&
\lesssim 
 \|a\|_{B^0_{\infty,\infty}}\left(\int_{\R}\left|\sum_j 2^{-j}|\Delta^{1/4} h^{j-4}|^2\right|^{q^{\prime}/2}\right)^{1/q^{\prime}}\left(\int_{\R}\left|\sum_j 2^{j}|b_j|^2\right|^{q/2}\right)^{1/q}\\[5mm]
 &&
 \lesssim 
  \|a\|_{B^0_{\infty,\infty}}\|b\|_{W^{1/2,q}}\,.
  \end{eqnarray*}
  
 {\bf $\bullet$ Estimate of $\|\Pi_3((\Delta^{1/4}a) b)\|_{L^q}\,.$}\par
 \begin{eqnarray*}
&& \|\sum_j\Delta^{1/4}a_j b_j\|_{L^q}  =\sup_{\|h\|_{L^{q^{\prime}}}\le 1}\int_{\R  h \sum_j }\Delta^{1/4} a_j b_j dx\\[5mm]
  &&
  =\sup_{\|h\|_{L^{q^{\prime}}}\le 1}\left[\int_{\R}  \sum_j \sum_{|k-j|\le 4 }  h_k (\Delta^{1/4}  a_j )b_j dx+\int_{\R}  \sum_j h^{j-4} (\Delta^{1/4}  a_j )b_j dx\right]
  \end{eqnarray*}
  We estimate the last term $\int_{\R}  \sum_j  
  h^{j-4} \Delta^{1/4}  a_j b_j dx$. \par
  To this purpose we show that $\sum_j\Delta^{1/4}(h^{j-4}b_j)\in h^1$ and the conclusion follows from
  the embedding $\dot H^{1/2}(\R)\hookrightarrow BMO(\R)\,.$
  We have
  \begin{eqnarray*}
  &&
  \int_{\R} \left( \sum_j 2^j |h^{j-4}  b_j |^2\right)^{1/2}dx\\[5mm]
  &&
  \lesssim
  \int_{\R} \sup_j |h^{j-4}|\left(\sum_j 2^j |b_j|^2\right)^{1/2} dx\\[5mm]
  &&
  \lesssim
  \left(\int_{\R}\sup_j|h^{j-4}|^{q^{\prime}}\right)^{1/q^{\prime}}
  \left(\int_{\R}(\sum_j 2^j |b_j|^2)^{q/2}\right)^{1/q}\\[5mm]
  &&
  \lesssim
  \|h\|_{L^{q^{\prime}}}\|b\|_{W^{1/2,q}}\,.
  \end{eqnarray*}

 {\bf $\bullet$ Estimate of $\|\Pi_2(\Delta^{1/4}(ab)-(\Delta^{1/4}a) b)\|_{L^q}\,.$}\par
 \begin{eqnarray*}
 &&\|\Pi_2(\Delta^{1/4}(ab)-(\Delta^{1/4}a) b)\|_{L^q}\\[5mm]
 && \simeq \sup_{\|h\|_{L^{q^{\prime}}}\le 1}\int_{\R}\sum_j h_j(\Delta^{1/4}(a_j b^{j-4})-\Delta a_j b^{j-4}) dx\\[5mm]
  && \simeq
  \sup_{\|h\|_{L^{q^{\prime}}}\le 1}\int_{\R}\sum_j b^{j-4}(\Delta^{1/4}(h_j) a_j-h_j\Delta a_j )dx\\
  &&
  \sup_{\|h\|_{L^{q^{\prime}}}\le 1}\int_{\R}\sum_j {\cal{F}}[b]^{j-4}(\eta)\int_{\R}{\cal{F}}[h]_j(\xi){\cal{F}}[a]_j(\eta-\xi)(|\xi|^{1/2}-|\eta-\xi|^{1/2}) d\xi\\[5mm]
  &&
  \cdots
  \\
  &&
  \lesssim
  \|a\|_{B^0_{\infty,\infty}}\|b\|_{W^{1/2,q}}\|h\|_{L^{q^{\prime}}}\,. \end{eqnarray*}\par
  This concludes the proof of Lemma \ref{lemmaprel4}\,.~\hfill$\Box$\par
  \medskip
  In the next Theorem we prove an estimate for the dual of the operator $F$ introduced in \rec{F}. It is defined as follows: given $Q\in L^2(\R)$,  $v\in \dot H^{1/2}(\R)$ we have
  $$F^{*}(Q,v)=\Delta^{1/4}(Q v)-\Delta^{1/4}{\mathcal{R}}({\mathcal{R}}(Q)v)\,.$$
  
\begin{Lemma}\label{Romega}
Let $Q\in L^2(\R)$,  $v\in \dot H^{1/2}(\R)$. Then
\begin{equation}
\|\Delta^{1/4}(Q v)-\Delta^{1/4}{\mathcal{R}}({\mathcal{R}}(Q)v)\|_{{\mathcal{H}}^{1}}\lesssim
\|Q\|_{L^2}\|v\|_{\dot H^{1/2}}\,.
\end{equation}
\end{Lemma}
{\bf Proof of Lemma \ref{Romega}\,.}\par
\noindent{\bf Estimate of  $\Pi_2(\Delta^{1/4}(Q,v))\,.$}
\begin{eqnarray}\label{P1Omega}
\|\Pi_1(\Delta^{1/4}(Q,v))\|_{{\cal{H}}^{1}}&=&\int_{\R}\left(\sum_{i=-\infty}^{+\infty} 2^i( Q^{i-4})^2(v_i)^2\right)^{1/2}dx\\[5mm]
&\lesssim &\int_{\R} |M(Q)|\left(\sum_{i=-\infty}^{+\infty} 2^i(v_i)^2\right)dx\nonumber\\[5mm]
&\lesssim& \|Q\|_{L^2}\|v\|_{\dot H^{1/2}}\,.\nonumber
\end{eqnarray}
The estimate of $\Pi_1(\Delta^{1/4}{\mathcal{R}}({\mathcal{R}}(Q)v))$ is analogous to
\rec{P1Omega}\,.
\par 
{\bf Estimate of  $\Pi_3(\Delta^{1/4}(Q,v))\,.$}
\begin{eqnarray}\label{P3Omega}
\|\Pi_3(Q,v)\|_{B^0_{1,1}}&\simeq&\sup_{\|h\|_{B^0_{\infty,\infty}}\le 1}\int_{R}
(Q_iv_i)\left[\Delta^{1/4} h^{i-6}+\sum_{t=h-5}^{i+6} \Delta^{1/4} h_t\right] dx\\
&
\lesssim& \sup_{\|h\|_{B^0_{\infty,\infty}}\le 1}\|h\|_{B^0_{\infty,\infty}}\int_{\R}2^{i/2}|Q_iv_i |dx\nonumber\\
&\lesssim&\left( \int_{\R}\sum_i 2^iv^2_idx\right)^{1/2}\left( \int_{\R}\sum_i Q^2_idx\right)^{1/2}=\|Q\|_{L^2}\|v\|_{\dot H^{1/2}}\,.\nonumber
 \end{eqnarray}
 The estimate of $\Pi_3(\Delta^{1/4}{\mathcal{R}}({\mathcal{R}}(Q)v))$ is analogous to
\rec{P3Omega}\,.
\par 
Finally one can easily check that
$$
\|\Pi_1(\Delta^{1/4}(Q v)-\Delta^{1/4}{\mathcal{R}}({\mathcal{R}}(Q)v))\|_{{\mathcal{H}}^{1}}=0\,.$$
This concludes the proof of Lemma \ref{Romega}\,.~\hfill$\Box$

  \end{document}